\documentclass[12pt]{article}
\usepackage{amsmath,amssymb,amsthm}
\usepackage{bbm}
\usepackage{comment}
\usepackage{hyperref}

\newcommand{\ds}{\displaystyle}
\newcommand{\fr}[2]{\frac{#1}{#2}}
\newcommand{\dfr}[2]{\dfrac{#1}{#2}}
\newcommand{\cd}{\cdot}
\newcommand{\cds}{\cdots}
\newcommand{\dsum}{\displaystyle \sum}
\newcommand{\ol}[1]{\overline{#1}}

\renewcommand{\l}{\left}
\renewcommand{\r}{\right}
\newcommand{\la}{\langle}
\newcommand{\ra}{\rangle}

\newcommand{\abs}[1]{\lvert{#1}\rvert}

\DeclareMathOperator*{\tensor}{\otimes}
\DeclareMathOperator*{\mystar}{\ast}
\DeclareMathOperator*{\mycirc}{\circ}

\newcommand{\Z}{\mathbb{Z}}
\newcommand{\C}{\mathbb{C}}
\newcommand{\R}{\mathbb{R}}
\newcommand{\N}{\mathbb{N}}

\newcommand{\res}{\mathrm{Res}}
\newcommand{\End}{\mathrm{End}}
\newcommand{\vir}{\mathrm{Vir}}
\newcommand{\ns}{\mathrm{NS}}
\newcommand{\ram}{\mathrm{R}}
\newcommand{\aut}{\mathrm{Aut}}
\newcommand{\wt}{\mathrm{wt}}

\newcommand{\sfr}[2]{\leavevmode\kern-.05em
  \raise.5ex\hbox{\the\scriptfont0 #1}\kern-.1em
  /\kern-.15em\lower.25ex\hbox{\the\scriptfont0 #2}\kern.02em}
\newcommand{\shf}{\sfr{1}{2}}

\newcommand{\Span}{\mathrm{Span}}
\newcommand{\Ker}{\mathrm{Ker}}
\newcommand{\w}{\omega}
\newcommand{\vac}{\mathbbm{1}}
\newcommand{\fusion}{\boxtimes}

\makeatletter
\@addtoreset{equation}{section}

\theoremstyle{plain}
\newtheorem{thm}{Theorem}[section]
\newtheorem{prop}[thm]{Proposition}
\newtheorem{lem}[thm]{Lemma}

\theoremstyle{definition}

\theoremstyle{remark}
\newtheorem{rem}[thm]{Remark}

\newcommand{\pf}{\noindent \textbf{Proof:}~~}

\title{3-dimensional Griess algebras and Miyamoto involutions}

\author{
   Ching Hung Lam\footnote{Partially supported by MoST grant 104-2115-M-001-004-MY3.}%
  \medskip\\
  {\small \it Institute of Mathematics, Academia Sinica, Taipei, Taiwan 10617}\\
 {\small \it  and National Center for Theoretical Sciences, Taipei, Taiwan}\\
  {\small e-mail: {\tt chlam@math.sinica.edu.tw}}
  \medskip\\
  Hiroshi Yamauchi\footnote{Partially supported by JSPS Grant-in-Aid for Young Scientists (B)
  No 24740027. }
  \medskip\\
  {\small \it Department of Mathematics,
  Tokyo Woman's Christian University}\\
  {\small \it 2-6-1 Zempukuji, Suginami-ku, Tokyo 167-8585, Japan}\\
  {\small e-mail: \texttt{yamauchi@lab.twcu.ac.jp}}
}
\date{}

\newcommand{\com}{\mathrm{Com}}
\renewcommand{\o}{\mathrm{o}}
\newcommand{\tw}{\mathrm{tw}}
\newcommand{\mydiv}{\,|\,}

\begin{document}

\maketitle

\baselineskip 6mm

\begin{abstract}
We consider a series of VOAs generated by 3-dimensional Griess algebras.
We will show that these VOAs can be characterized by their 3-dimensional Griess 
algebras and their structures are uniquely determined.
As an application, we will determine the groups generated by the Miyamoto involutions associated to Virasoro vectors 
of our VOAs.
\end{abstract}

\tableofcontents

\section{Introduction}

Let $L(c_n^0,0)$ be the unitary simple Virasoro vertex operator algebra (VOA) 
of central charge $c_n^0$, where 
\[ 
  c_n^0=1-\dfr{6}{(n+2)(n+3)},~~n=1,2,3,\dots.
\]
Let $V$ be  a VOA containing a sub VOA $U$ isomorphic to $L(c_n^0,0)$ 
and let $e$ be the Virasoro vector of $U$.
It was shown by Miyamoto \cite{Mi1} that one can define an automorphism of $V$ 
based on the fusion rules of $L(c_n^0,0)$-modules.
Namely, the zero-mode $\o(e)$ of $e$ acts on $V$ semisimply and the linear map
$\tau_e=(-1)^{4(n+2)(n+3)\o(e)}$ is well-defined on $V$. The map $\tau_e$ gives rise to an 
automorphism of $V$ and we will call it the \textit{Miyamoto involution} associated to $e$.
When $\tau_e$ is trivial,  one can define another Miyamoto 
involution $\sigma_e$ on a certain sub VOA of $V$ using a similar method.
(See Theorem \ref{thm:2.3} for detail.)
When the involution $\sigma_e$ is well-defined on the whole space $V$, 
such a Virasoro vector $e$ is said to be of \textit{$\sigma$-type} on $V$.
Miyamoto involutions appear naturally as symmetries of some sporadic 
finite simple groups in the VOA theory.

The most important and interesting case would be the first member $c_1^0=1/2$ 
of the unitary series.
It is shown in \cite{Mi1} that there exists a one-to-one correspondence 
between the 2A-elements of the Monster simple group and sub VOAs of the 
moonshine VOA $V^\natural$ isomorphic to $L(\shf,0)$ via the Miyamoto involutions.
However, the case $c=c_1^0$ is not the only important case.
It is shown in \cite{HLY1} that there exists a one-to-one correspondence 
between the 2A-elements of the Baby Monster simple group and 
$c=c_2^0$ Virasoro vectors of $\sigma$-type in the VOA $V\!\!B^\natural$. 
An injective correspondence between the 2C-elements of the largest Fischer 
3-transposition group and $c=c_4^0$ Virasoro vectors of $\sigma$-type in the VOA 
$V\!\!F^\natural$ is also shown in \cite{HLY2}.
Recently, in \cite{LY}, the authors found new correspondences between the transpositions 
of the second and the third largest Fischer 3-transposition groups 
and $c=c_5^0$ Virasoro vectors and $c=c_6^0$ Virasoro vectors of $\sigma$-type, respectively.
Those relations between involutions and Virasoro vectors are referred to as 
the Conway-Miyamoto correspondence in (loc.~cit.).

The purpose of this paper is to study the relations between Miyamoto involutions 
associated to two mutually orthogonal Virasoro vectors of $c=c_n^0$ and $c=c_{n+1}^0$.
Straightly speaking, there will be no special relations between their Miyamoto involutions 
if the Virasoro vectors are just mutually orthogonal. 
Therefore, we need to add extra assumptions and consider their extensions. 
Namely, we will consider a VOA $A(\shf,c_n^1)$ generated by a 3-dimensional 
Griess algebra which is spanned by two mutually orthogonal vectors and one common 
highest weight vector.
This VOA has already been  considered in \cite{A2,LS} and occurs naturally 
in the moonshine VOA $V^\natural$ \cite{DLMN}.
In this paper, we will give a characterization of this VOA based on 
the 3-dimensional Griess algebra in Theorem 3.5.
In the description, we will use commutant superalgebras studied in \cite{Y1} 
and the $N=1$ super Virasoro algebras. 
A main observation is that 
\[
  c_n^0+c_{n+1}^0-\frac{1}2 
  =  c_n^1=\dfr{3}{2}\l(1-\dfr{8}{(n+2)(n+4)}\r),~~ \text{ for } n=1,2,3,\dots,
\]
while the number $c_n^1$ is  the central charge of the unitary $N=1$ super Virasoro algebra.
(The superscript in $c_n^1$ means ``$N=1$'', whereas that in $c_n^0$ means 
``$N=0$'', i.e., non-super case.)
The VOA $A(\shf,c_n^1)$ has three Virasoro frames and we can consider 
associated Miyamoto involutions.
Two of them consist of mutually orthogonal $c=c_n^0$ and $c=c_{n+1}^0$ 
Virasoro vectors.
We will classify all the irreducible modules over $A(\shf,c_n^1)$ and 
describe their decompositions with respect to these three Virasoro frames.
As an application, we will determine the group generated by Miyamoto 
involutions associated to these Virasoro frames of $A(\shf,c_n^1)$ in Theorem 4.6.
In Theorem 4.8, we will present certain inductive relations 
between Miyamoto involutions associated to $c=c_n^0$ and $c=c_{n+1}^0$ Virasoro 
vectors of $A(\shf,c_n^1)$ when $n$ is odd.
It is worthy to mention that our characterization is very simple and 
easy to check in practice.
Actually, this work is a part of study of Fischer 3-transposition groups 
in \cite{LY} and the results in this paper will be crucial in the discussion 
in Section 5 of (loc.~cit.).

The organization of this article is as follows. 
We first review some basic notation and terminology about VOAs. 
In Section 2, we review some basic facts and results for the unitary series of 
Virasoro VOAs and the unitary series of the $N=1$ super Virasoro algebras. 
In Section 3, we study VOAs generated by their $3$-dimensional Griess algebras. 
In Theorem \ref{thm:3.5}, we will gives a characterization of 
such a VOA using its Griess algebra. 
In Section 4, we introduce a series of VOAs $A(\shf,c_n^1)$ which satisfy 
Theorem \ref{thm:3.5}.
We study three Virasoro frames of $A(\shf,c_n^1)$ and classify
all its irreducible modules and their decompositions with respect to 
three Virasoro frames.
We will then determine the actions of Miyamoto 
involutions on VOAs which contain $A(\shf,c_n^1)$ as a sub VOA in Theorems \ref{thm:4.6} and \ref{thm:4.8}.
This is the main result of this paper.
In Appendix, we will prove the $\Z_2$-rationality of the unitary series of 
the Ramond algebra, one of the $N=1$ super Virasoro algebras.

\paragraph{Acknowledgement.}
Part of this work has been done while the authors were staying 
at National Taitung University in March 2015,  
at Mathematisches Forschungsinstitut Oberwolfach in March 2015, 
and at Sichuan University in September 2015.
They gratefully acknowledge the hospitality there.
H.Y. thanks Dra\v{z}en Adamovi\'{c} and Yoshiyuki Koga for valuable comments and references.

\paragraph{Notation and terminology.}
In this paper, we will work mainly over the complex number field $\C$ but 
sometimes we also consider real forms.
If $V$ is a superspace,  we denote its $\Z_2$-grading by $V=V^{[0]}\oplus V^{[1]}$.
A VOA $V$ is called \textit{of OZ-type} if it has the grading 
$V=\oplus_{n\geq 0}V_n$ such that $V_0=\C \vac$ and $V_1=0$.
We will mainly consider VOAs of OZ-type.
In this case, $V$ has a unique invariant bilinear form such that $(\vac\mydiv \vac)=1$.
A real form $V_\R$ of $V$ is called \emph{compact} if the associated bilinear form 
is positive definite.
For a subset $A$ of $V$, the subalgebra generated by $A$ is denoted by $\la A\ra$.
For $a\in V_n$, we define $\wt(a)=n$.
We expand $Y(a,z)=\sum_{n\in \Z}a_{(n)}z^{-n-1}$ for $a\in V$ and define 
its \textit{zero-mode} by $\o(a):=a_{(\wt(a)-1)}$ if $a$ is homogeneous and extend linearly.
The weight two subspace $V_2$ of a VOA $V$ of OZ-type carries a structure of a commutative 
algebra equipped with the product $\o(a)b=a_{(1)}b$ for $a$, $b\in V_2$.
This algebra is called the \emph{Griess algebra} of $V$.
A \emph{Virasoro vector} is a vector $a\in V_2$ such that $a_{(1)}a=2a$.
In this case the subalgebra $\la a\ra$ is isomorphic to a Virasoro VOA with 
the central charge $c=2(a\mydiv a)$.
A Virasoro vector $e\in V$ is called \textit{simple} if it generates a simple 
Virasoro sub VOA. Notice that $e$ is always simple if it is taken from a compact real form of $V$.
A simple $c=1/2$ Virasoro vector is called an \emph{Ising vector}.
A Virasoro vector $\w$ is called the \emph{conformal vector} of $V$ if each graded 
subspace $V_n$ agrees with $\Ker_V\, (\o(\w)-n)$ and satisfies 
$\w_{(1)}a =a_{(-2)}\vac$ for all $a\in V$.
When $V$ is of OZ-type,  the half of the conformal vector gives the unit of the Griess algebra and hence it is uniquely determined.
We write $L(n)=\w_{(n+1)}$ for $n\in \Z$.
A \emph{sub VOA $(W,e)$} of $V$ is a pair of a subalgebra $W$ of $V$ and
a Virasoro vector $e\in W$ such that $e$ is the conformal vector of $W$.
Usually we omit to denote $e$ and simply call $W$ a sub VOA.
A sub VOA $W$ of $V$ is said to be \textit{full} if $V$ and $W$ shares the same conformal vector.
The \textit{commutant subalgebra} of a sub VOA $(W,e)$ of $V$ is defined by 
$\com_V W=\Ker_V e_{(0)}$ (cf.~\cite{FZ}).
For an automorphism $\sigma$ of $V$ and a $V$-module $(M,Y_M(\cd,z))$, we define its 
$\sigma$-conjugate by $(M,Y_M^\sigma(\cd,z))$ where $Y_M^\sigma(a,z)=Y_M(\sigma a,z)$ for $a\in V$.
A module $M$ is called \textit{$\sigma$-invariant} or \textit{$\sigma$-stable} if 
$M$ is isomorphic to its $\sigma$-conjugate.

\section{Virasoro VOAs and SVOAs}

\subsection{Unitary series of the Virasoro algebra}

Let $\vir=\oplus_{n\in \Z} \C L(n)\oplus \C z$ be the Virasoro algebra. 
The irreducible highest weight module over $\vir$ with central charge $c$ and 
highest weight $h$ will be denoted by $L(c,h)$. 
Let
\begin{equation}\label{eq:2.2}
\begin{array}{l}
  c_n^0:= 1-\dfr{6}{(n+2)(n+3)},~~~ n=1,2,3,\dots,
  \medskip\\
  h_{r,s}^{(n)}:=\dfr{(r(n+3)-s(n+2))^2-1}{4(n+2)(n+3)},~~ 1\leq r\leq n+1,~~ 1\leq s\leq n+2.
\end{array}
\end{equation}
Then $L(c_n^0,0)$ is a rational $C_2$-cofinite VOA and $L(c_n^0,h_{r,s}^{(n)})$, 
$1\leq s\leq r\leq n+1$, exhaust the set of inequivalent irreducible $L(c_n^0,0)$-modules  
(cf.~\cite{DMZ,W}).
Note that $h^{(n)}_{r,s}=h^{(n)}_{n+2-r,n+3-s}$.
It is known that all irreducible $L(c_n^0,0)$-modules have compact real forms and 
$L(c_n^0,h_{r,s}^{(n)})$ are usually called the unitary series of the Virasoro algebra.
The fusion rules are also known and given as follows.
\begin{equation}\label{eq:2.3}
  L(c_n^0,h_{r,s}^{(n)})\fusion L(c_n^0,h_{r',s'}^{(n)})
  = \dsum_{1\leq i \leq M \atop 1\leq j \leq N}
    L(c_n^0,h_{\abs{r-r'}+2i-1,\abs{s-s'}+2j-1}^{(n)}), 
\end{equation}
where $M=\min\{ r,r',n+2-r,n+2-r'\}$ and $N=\min\{ s,s',n+3-s,n+3-s'\}$.

Let $V$ be a VOA $V$ and $e$ a Virasoro vector of $V$.
Suppose $e$ generates a simple $c=c_n^0$ Virasoro sub VOA of $V$.
Then $V$ is an $\la e\ra$-module and one has an isotypical decomposition 
\begin{equation}\label{eq:2.4}
  V=\bigoplus_{1\leq s\leq r\leq n+1} V_e[h^{(n)}_{r,s}],
\end{equation}
where $V_e[h]$ is a sum of all irreducible $\la e\ra$-submodules isomorphic to $L(c_n^0,h)$.
Moreover, the zero-mode $\o(e)=e_{(1)}$ of $e$ acts semisimply on $V$ .

\begin{thm}[\cite{Mi1}]\label{thm:2.1}
  Let $e$ be a simple $c=c_n^0$ Virasoro vector of $V$.
  Then the linear map $\tau_e=(-1)^{4(n+2)(n+3)\o(e)}$ defines an automorphism in $\aut(V)$.
  By \eqref{eq:2.2},  $\tau_e$ acts on $V_e[h_{r,s}^{(n)}]$ as $(-1)^{r+1}$ 
if $n$ is even and as $(-1)^{s+1}$ if $n$ is odd.
\end{thm}

Set
\begin{equation}\label{eq:2.5}
  P_n:=\begin{cases}
    \{ h_{1,s}^{(n)} \mid 1\leq s\leq n+2\} & \mbox{if $n$ is even,}
    \medskip\\
    \{ h_{r,1}^{(n)} \mid 1\leq r\leq n+1\} & \mbox{if $n$ is odd.}
  \end{cases}
\end{equation}
It follows from the fusion rules in \eqref{eq:2.3} that the subspace 
$V_e[P_n]=\bigoplus_{h\in P_n} V_e[h]$
forms a subalgebra of $V$.
We say $e$ is \emph{of $\sigma$-type on $V$} if $V=V_e[P_n]$.
Again by the fusion rules in \eqref{eq:2.3}, we have the following $\Z_2$-symmetry.

\begin{thm}[\cite{Mi1}]\label{thm:2.3}
  The linear map 
  \[
    \sigma_e :=\begin{cases}
      (-1)^{s+1} ~~\mbox{on~ $V_e[h_{1,s}^{(n)}]$~ if~ $n$~ is even,}
      \medskip\\
      (-1)^{r+1} ~~\mbox{on~ $V_e[h_{r,1}^{(n)}]$~ if~ $n$~ is odd,} 
    \end{cases}
  \]
  defines an automorphism of $V_e[P_n]$. 
\end{thm}

\subsection{Unitary series of the $N=1$ super Virasoro algebra}\label{sec:2.2}

There are two extensions of the Virasoro algebra to Lie superalgebras 
called the $N=1$ super Virasoro algebras:   
$\ns=\ns^{[0]}\oplus \ns^{[1]}$ and $\ram=\ram^{[0]}\oplus \ram^{[1]}$, 
where 
\[
  \ns^{[0]}=\ram^{[0]}=\vir=\bigoplus_{n\in \Z} \C L(n)\oplus \C z,~~~
  \ns^{[1]}=\bigoplus_{r\in \Z+1/2} \C G(r), ~~~~
  \ram^{[1]}=\bigoplus_{r\in \Z} \C G(r).
\]
In addition to the Virasoro algebra relations, they also satisfy the following relations:
\begin{equation}\label{eq:2.7}
\begin{array}{l}
  [L(m),G(r)]=\l(\dfr{1}{2}m-r\r) G(m+r), ~~~[z,G(r)]=0, \text{ and } 
  \medskip\\
  {}[G(r),G(s)]=2L(r+s)+\delta_{r+s,0}\l( r^2-\dfr{1}{4}\r) \dfr{z}{3}.
\end{array}
\end{equation}
$\ns$ is called the \textit{Neveu-Schwarz algebra} and $\ram$ is called 
the \textit{Ramond algebra}.

\paragraph{Neveu-Schwarz sectors.}
The Neveu--Schwarz algebra has a standard triangular decomposition 
$\ns=\ns_-\oplus \ns_0\oplus \ns_+$ such that 
\[
  \ns_{\pm}=\bigoplus_{\pm n>0} \C L(n) \oplus \bigoplus_{\pm r>0} \C G(r),~~
  \ns_0=\C L(0) \oplus \C z.
\]
For $c, h\in \C$, let $\C v_{c,h}$ be a one-dimensional module over $\ns_0\oplus \ns_+$ 
defined by 
\[
  zv_{c,h}=cv_{c,h},~~
  L(0)v_{c,h}=h v_{c,h},~~
  \ns_+v_{c,h}=0,
\]
and define the Verma module over the Neveu-Schwarz algebra $M_\ns(c,h)$ with central charge 
$c$ and highest weight $h$ by the induced module.
We choose the $\Z_2$-grading of $M_\ns(c,h)=M_\ns(c,h)^{[0]}\oplus M_\ns(c,h)^{[1]}$ so that  
$v_{c,h}\in M_\ns(c,h)^{[0]}$.
We denote the unique simple quotient of $M_\ns(c,h)$ by $L_\ns(c,h)$ which is called 
the NS-sector.
Set
\begin{equation}\label{eq:2.8}
  L(z)=\dsum_{n\in \Z} L(n)z^{-n-2}, \quad \text{ and } \quad 
  G_\ns(z)=\dsum_{r\in \Z+1/2} G(r)z^{-r-3/2}.
\end{equation}
Then $L(z)$ and $G_\ns(z)$ are local fields on $M_\ns(c,h)$.
Namely, they are elements of $\End(M_\ns(c,h))[\![z,z^{-1}]\!]$ and  one has the following OPEs:
\begin{equation}\label{eq:2.9}
\begin{array}{l}
  L(z)L(w)\sim \dfr{c}{2(z-w)^4}+\dfr{2L(w)}{(z-w)^2}+\dfr{\partial_w L(w)}{z-w},
  \medskip\\
  L(z)G_\ns(w)\sim \dfr{3G_\ns(w)}{2(z-w)^2}+\dfr{\partial G_\ns(w)}{z-w},
  \medskip\\
  G_\ns(z)G_\ns(w)\sim \dfr{2c}{3(z-w)^3}+\dfr{2L(w)}{z-w} .
\end{array}
\end{equation}
We also have the derivation relations
\[
[L(-1),L(z)]=\partial_z L(z)\quad \text{ and }\quad  [L(-1),G_\ns(z)]=\partial_z G_\ns(z).
\]
Therefore, $L(z)$ and $G_\ns(z)$ generate a vertex superalgebra inside 
$\End(M_\ns(c,h))[\![z,z^{-1}]\!]$.
Set $\ol{M}_\ns(c,h):=M_\ns(c,0)/\la L(-1)v_{c,0}\ra$.
We denote the images of $v_{c,0}$, $L(-2)v_{c,0}$ and $G(-\sfr{3}{2})v_{c,0}$ in the quotient 
$\ol{M}_\ns(c,0)$ by $\vac$, $\w$ and $\tau$, respectively.
Then $\ol{M}_\ns(c,0)$ carries a unique structure of a vertex operator superalgebra 
such that $\vac$ is the vacuum vector, 
$\w=L(-2)\vac=\fr{1}{2}G(-\sfr{1}{2})G(-\sfr{3}{2})\vac$ is the conformal vector 
satisfying $Y(\w,z)=L(z)$ and $Y(\tau,z)=G_\ns(z)$ by Theorem 4.5 of \cite{K}.
Note that $\la G(-\sfr{1}{2})v_{c,0}\ra=\la L(-1)v_{c,0}\ra$ in $M_\ns(c,0)$ as 
$G(-\sfr{1}{2})^2=L(-1)$ and $[G(\sfr{1}{2}),L(-1)]=G(-\sfr{1}{2})$.
It is clear that $M_\ns(c,h)$ is a $\ol{M}_\ns(c,0)$-modules such that 
$Y(\tau,z)=G_\ns(z)$ on $M_\ns(c,h)$ and its simple quotient $L_\ns(c,h)$ is an irreducible 
$\ol{M}_\ns(c,0)$-module.
In particular, the simple quotient $L_\ns(c,0)$ is a simple SVOA.

\paragraph{Ramond sectors.}
The Ramond algebra also has a standard triangular decomposition 
$\ram=\ram_-\oplus \ram_0\oplus \ram_+$ such that 
\[
  \ram_{\pm}=\bigoplus_{\pm n>0} \C L(n) \oplus \bigoplus_{\pm r>0} \C G(r),~~
  \ram_0=\C L(0) \oplus \C z \oplus \C G(0).
\]
Since $G(0)^2=L(0)-z/24$, the subalgebra $\ram_0$ is not $\Z_2$-homogeneous in this case.
For $c$, $d\in \C$,  let $\C v_{c,d}$ be a one-dimensional module over $\ram_0\oplus \ram_+$ 
defined by 
\[
  zv_{c,d}=cv_{c,d},~~
  G(0)v_{c,d}=d v_{c,d},~~
  \ram_+v_{c,d}=0,
\]
and define the Verma module over the Ramond algebra $M_\ram(c,d)$ with central 
charge $c$ and top weight $d$ by the induced module.
Note that $L(0)v_{c,d}=(d^2+c/24)v_{c,d}$.
We call the eigenvalue $d^2+c/24$ the \textit{highest weight} of $M_\ram(c,d)$.
Our notion of Verma modules is unusual in the sense that $G(0)$ always acts 
semisimply on highest weight vectors and there is no canonical superspace 
structure.
We denote the unique simple quotient of $M_\ram(c,d)$ by $L_\ram(c,d)$ which is 
called the Ramond sector.
Set
\begin{equation}\label{eq:2.10}
  L(z)=\dsum_{n\in \Z} L(n) z^{-n-2}, \quad \text{ and } \quad
  G_\ram(z)=\dsum_{n\in \Z} G(n)z^{-n-3/2}.
\end{equation}
Then $L(z)$ and $G_\ram(z)$ are local $\Z_2$-twisted fields on $M_\ram(c,d)$. 
Namely, one has the same OPEs as in \eqref{eq:2.9}.
We also have the derivation relations $[L(-1),L(z)]=\partial_z L(z)$ and 
$[L(-1),G_\ram(z)]=\partial_z G_\ram(z)$.
Therefore, $L(z)$ and $G_\ram(z)$ generate a vertex superalgebra 
inside $\End(M_\ram)[\![z^{1/2},z^{-1/2}]\!]$ (cf.~\cite{Li1}).
Let $\theta=(-1)^{2L(0)}$ be the canonical $\Z_2$-symmetry (supersymmetry) of $\ol{M}_\ns(c,0)$.
One can directly verify that $M_\ram(c,d)$ is a $\theta$-twisted 
$\ol{M}_\ns(c,0)$-module and its simple quotients $L_\ram(c,d)$ 
is an irreducible $\theta$-twisted $\ol{M}_\ns(c,0)$-module.
It follows that $L_\ram(c, d)$ and $L_\ram(c, -d)$ are $\theta$-conjugate to each other.

\paragraph{Unitary series.}
For $n=1,2,3,\dots$, let
\begin{equation}\label{eq:2.11}
\begin{array}{l}
  c_n^1=\dfr{3}{2}\l( 1-\dfr{8}{(n+2)(n+4)}\r),
  \medskip\\
  \Delta_{r,s}^{(n)}=\dfr{s(n+2)-r(n+4)}{\sqrt{8(n+2)(n+4)}},~~~~~
  \medskip\\
  h_{r,s,p}^{(n)}=\big(\Delta_{r,s}^{(n)}\big)^2+\dfr{c_n^1}{24}-\dfr{1}{16}\delta_{p,0}
  =\dfr{(r(n+4)-s(n+2))^2-4}{8(n+2)(n+4)}+\dfr{p}{8},
\end{array}
\end{equation}
where $1\leq r\leq n+1$, $1\leq s\leq n+3$, $p=0$ or $1/2$ and $r-s\equiv 2p \mod 2$.
Note that $\Delta_{n+2-r,n+4-s}^{(n)}=-\Delta_{r,s}^{(n)}$ and $h_{n+2-r,n+4-s,p}^{(n)}=h_{r,s,p}^{(n)}$.

\begin{thm}[\cite{A1,IK,KW,M}]\label{thm:2.4}
  The $N=1$ Virasoro SVOA $L_{\mathrm{NS}}(c_n^1,0)$ is rational and $\Z_2$-rational.
  The irreducible representations are as follows.
  \\
  (1)~ The NS-sectors $L_{\mathrm{NS}}(c_n^1,h_{r,s,0}^{(n)})$, $1\leq r\leq m+1$, 
  $1\leq s\leq m+3$, $r\equiv s\mod 2$, are all the irreducible untwisted 
  $L_{\mathrm{NS}}(c_n^1,0)$-modules.
  \\
  (2)~ The R-sectors $L_{\ram}(c_n^1,\Delta_{r,s}^{(n)})$, $1\leq r\leq m+1$, 
  $1\leq s\leq m+3$, $r\equiv s+1\mod 2$, are all the irreducible $\Z_2$-twisted 
  $L_{\mathrm{NS}}(c_n^1,0)$-modules.
\end{thm}

\pf
The irreducible untwisted $L_{\mathrm{NS}}(c_n^1,0)$-modules are classified 
in \cite{KW,A1} and the rationality is established in \cite{A1}.
The irreducible $\Z_2$-twisted $L_{\mathrm{NS}}(c_n^1,0)$-modules are classified in \cite{M}.
The $\Z_2$-rationality will be given in Appendix.
\qed

\begin{rem}
The NS-sectors $L_\ns(c_n^1,h_{r,s,0}^{(n)})$ and 
the R-sectors $L_\ram(c_n^1,\Delta_{r,s}^{(n)})$ have compact real forms and 
are called the \textit{unitary series} of the $N=1$ super Virasoro algebra.
\end{rem}

\begin{rem}
If $n$ is even, the top weight $\Delta_{(n+2)/2,(n+4)/2}^{(n)}=0$ is the fixed point 
of the $\Z_2$-symmetry $\Delta_{n+2-r,n+4-s,}^{(n)}=-\Delta_{r,s}^{(n)}$ and the corresponding 
representation $L_\ram(c_n^1,0)$ is $\theta$-stable. 
Therefore, among the irreducible R-sectors in (2) of Theorem \ref{thm:2.4}, the representation 
$L_\ram(c_n^1,0)$ has a distinguished property such that it is not irreducible as an 
$L_\ns(c_n^1,0)^{[0]}$-module while others are still irreducible over $L_\ns(c_n^1,0)^{[0]}$.
\end{rem}

\section{VOAs with 3-dimensional Griess algebras}

\subsection{Free Majorana fermion and $c=1/2$ Virasoro SVOA}\label{sec:3.1}

Let $C$ be the associative algebra generated by $\psi_r$, $r\in \Z+1/2$, subject to 
the relations
\begin{equation}\label{eq:3.1}
  \psi_r\psi_s+\psi_s\psi_r=\delta_{r+s,0},~~~r,s\in \Z+1/2.
\end{equation}
The Fock representation $F$ of $C$ is a cyclic $C$-module generated by 
$\vac$ with relations $\psi_r \vac=0$ for $r>0$.
Then $F$ has a natural $\Z_2$-grading $F=F^{[0]}\oplus F^{[1]}$ with 
\begin{equation}\label{eq:3.2}
  F^{[i]}=\Span_\C\{ \psi_{-r_1}\cds \psi_{-r_k}\vac \mid r_1>\cds >r_k>0,~~ k \equiv i \mod 2\} .
\end{equation}
The generating series $ \psi(z):=\sum_{r\in \Z+1/2}\psi_r z^{-r-1/2}$ is an odd field on $F$ 
and satisfies the following locality 
\begin{equation}\label{eq:3.3}
  (z_1-z_2) [\psi(z_1),\psi(z_2)]_+=0.
\end{equation}
Therefore, $\psi(z)$ generates a vertex superalgebra inside $\End(F)[\![z,z^{-1}]\!]$ and 
$F$ can be equipped with a unique structure of a vertex superalgebra 
such that $\vac$ is the vacuum vector and $Y(\psi_{-1/2}\vac,z)=\psi(z)$ (cf.~\cite{K}).
The vector $\w=\fr{1}{2}\psi_{-3/2}\psi_{-1/2}\vac$ provides the conformal vector of 
central charge 1/2 and we have the isomorphisms
\begin{equation}\label{eq:3.4}
  F^{[0]}\cong L(\shf,0),~~~
  F^{[1]}\cong L(\shf,\shf)
\end{equation}
as $\la \w\ra$-modules (cf.~\cite{KR}).

Let $C_{\mathrm{tw}}$ be the associative algebra generated by $\phi_n$, $n\in \Z$, 
subject to the relations
\begin{equation}\label{eq:3.5}
  \phi_m\phi_n +\phi_n\phi_n=\delta_{m+n,0},~~~m,n\in \Z .
\end{equation}
The Fock representation $F_{\mathrm{tw}}$ of $C_{\mathrm{tw}}$ is a cyclic $C_{\mathrm{tw}}$-module 
generated by $v_{1/16}$ with relations $\phi_n v_{1/16}=0$ for $n>0$.
Set $F_{\mathrm{tw}}^\pm:=\la v_{1/16}\pm \sqrt{2}\,\phi_0v_{1/16}\ra$.
Then we have a decomposition $F_{\mathrm{tw}}=F_{\mathrm{tw}}^+\oplus F_{\mathrm{tw}}^-$ as 
a $C_{\mathrm{tw}}$-module.
The generating series $\phi(z):=\sum_{n\in \Z}\phi_n z^{-n-1/2}$ is an odd $\Z_2$-twisted 
field on $F_{\mathrm{tw}}$ and satisfies the following locality 
\begin{equation}\label{eq:3.6}
  (z_1-z_2) [\phi(z_1),\phi(z_2)]_+=0.
\end{equation}
Therefore, $\phi(z)$ generates a vertex superalgebra inside 
$\End(F_{\mathrm{tw}})[\![z^{1/2},z^{-1/2}]\!]$ 
and $F_{\mathrm{tw}}$ can be equipped with a unique structure of a $(-1)^{2L(0)}$-twisted 
$F$-module such that $Y(\psi_{-1/2}\vac,z)=\phi(z)$ (cf.~\cite{Li1}).
As $F^{[0]}$-modules, we have the isomorphisms (cf.~\cite{KR})
\begin{equation}\label{eq:3.7}
  F_{\mathrm{tw}}^+\cong F_{\mathrm{tw}}^-\cong L(\shf,\sfr{1}{16}), 
\end{equation}
whereas $F_{\mathrm{tw}}^+$ and $F_{\mathrm{tw}}^-$ are inequivalent $(-1)^{2L(0)}$-twisted 
$F$-modules since the zero-mode $o(\psi_{-1/2}\vac)=\phi_0$ acts as $\pm 2^{-1/2}$ 
on the top levels of $F_{\mathrm{tw}}^\pm$.
Indeed, $F_{\mathrm{tw}}^\pm$ are mutually $(-1)^{2L(0)}$-conjugate $F$-modules.

Let $V=V^{[0]}\oplus V^{[1]}$ be an SVOA with a non-trivial odd part.
A tensor product $F\tensor V$ is an SVOA with $\Z_2$-grading
\begin{equation}\label{eq:3.8}
  (F\tensor V)^{[0]}=F^{[0]}\tensor V^{[0]} \oplus F^{[1]}\tensor V^{[1]},~~~
  (F\tensor V)^{[1]}=F^{[0]}\tensor V^{[1]} \oplus F^{[1]}\tensor V^{[0]}.
\end{equation}
Let $M=M^{[0]}\oplus M^{[1]}$ be an untwisted $V$-module.
Then $F\tensor M$ is an $F\tensor V$-module and its $\Z_2$-homogeneous parts
\begin{equation}\label{eq:3.9}
  (F\tensor M)^{[0]}=F^{[0]}\tensor M^{[0]}\oplus F^{[1]}\tensor M^{[1]},~~~
  (F\tensor M)^{[1]}=F^{[0]}\tensor M^{[1]}\oplus F^{[1]}\tensor M^{[0]},
\end{equation}
are $(F\tensor V)^{[0]}$-submodules.
Let $N$ be a $\Z_2$-twisted $V$-module.
Then tensor products 
\begin{equation}\label{eq:3.10}
  F_{\mathrm{tw}}^\pm\tensor N
\end{equation}
are untwisted $(F\tensor V)^{[0]}$-modules.
In this way, given an untwisted or $\Z_2$-twisted $V$-module, we can construct a 
$(F\tensor V)^{[0]}$-module.
In the next subsection, we will show that there is a canonical reverse construction 
of $V$-modules from $(F\tensor V)^{[0]}$-modules.

\subsection{Commutant SVOAs}

Let $V$ be a VOA and $e$ an Ising vector of $V$.
Let $M$ be a $V$-module.
For $h=0$, $1/2$, $1/16$, we set
\begin{equation}\label{eq:3.11}
  T_{e,M}(h):=\{ a\in V \mid e_{(1)}a=ha\} .
\end{equation}
Then we have the isotypical decomposition 
\begin{equation}\label{eq:3.12}
  M = L(\shf,0)\tensor T_{e,M}(0)
  \oplus L(\shf,\shf)\tensor T_{e,M}(\shf)
  \oplus L(\shf,\sfr{1}{16})\tensor T_{e,M}(\sfr{1}{16})
\end{equation}
of $M$ as $\la e\ra$-module. 
Since $\la e\ra\cong L(\shf,0)$ is rational, its zero-mode $\o(e)$ acts on $M$ semisimply. 
Therefore, the Miyamoto involution $\tau_e=(-1)^{48\o(e)}$ is also well-defined on $M$ and 
$M$ is $\tau_e$-stable. 
Set
\begin{equation}\label{eq:3.13}
\begin{array}{l}
  M^{\la \tau_e\ra} = L(\shf,0)\tensor T_{e,M}(0)
  \oplus L(\shf,\shf)\tensor T_{e,M}(\shf),
  \medskip\\
  M^{\la -\tau_e\ra}= L(\shf,\sfr{1}{16})\tensor T_{e,M}(\sfr{1}{16}).
\end{array}
\end{equation}
Then $M^{\la \pm \tau_e\ra}$ are $V^{\la \tau_e\ra}$-submodules.
It is known that $T_{e,V}(0)$ is the commutant subalgebra of $\la e\ra$ in $V$ and 
$T_{e,M}(h)$, $h=1/2$, $1/16$, are $T_{e,V}(0)$-modules (cf.~\cite{FZ,Y1}).
As $L(\shf,0)$ can be extended to an SVOA $L(\shf,0)\oplus L(\shf,\shf)$, 
the commutant $T_{e,V}(0)$ can also be extended to an SVOA.

\begin{prop}[Theorem 2.2 of \cite{Y1}]\label{prop:3.1}
  Let $V$ be a VOA and $e$ an Ising vector of $V$ and suppose $T_{e,V}(\shf)\ne 0$.
  Then there exists an SVOA structure on $T_{e,V}(0)\oplus T_{e,V}(\shf)$ 
  which is an extension of the commutant sub VOA $T_{e,V}(0)$ such that the even part 
  \[
    \l(\l( L(\shf,0)\oplus L(\shf,\shf)\r)\tensor 
    \l( T_{e,V}(0)\oplus T_{e,V}(\shf)\r)\r)^{[0]}
  \]
  of a tensor product of SVOAs is isomorphic to the sub VOA $V^{\la \tau_e\ra}$.
  If $V$ is simple, then $T_{e,V}(0)\oplus T_{e,V}(\shf)$ is a simple SVOA and 
  $T_{e,V}(0)$ is a simple sub VOA.
\end{prop}

Next theorem provides a construction of $T_{e,V}(0)\oplus T_{e,V}(\shf)$-modules 
from $V^{\la \tau_e\ra}$-modules, which is a sort of reverse of \eqref{eq:3.9} and \eqref{eq:3.10}.

\begin{prop}\label{prop:3.2}
  Let $V$ be a VOA and $e$ an Ising vector of $V$ such that $T_e(\shf)\ne 0$.
  Let $M$ be a $V^{\la \tau_e\ra}$-module.
  Decompose $M$ as in \eqref{eq:3.14} and define $M^{\la \pm \tau_e\ra}$ as in \eqref{eq:3.15}.
  \\
  (1) The space $T_{e,M}(0)\oplus T_{e,M}(\shf)$ forms an  untwisted 
  $T_{e,V}(0)\oplus T_{e,V}(\shf)$-module such that $M^{\la \tau_e\ra}$ as a $V^{\la \tau_e\ra}$-module 
  is isomorphic to one of the $\Z_2$-homogeneous parts of the tensor product of the adjoint 
  module of $L(\shf,0)\oplus L(\shf,\shf)$ and the $T_{e,V}(0)\oplus T_{e,V}(\shf)$-module 
  $T_{e,M}(0)\oplus T_{e,M}(\shf)$.
  If $M^{\la \tau_e\ra}$ is an irreducible $V^{\la \tau_e\ra}$-module, then 
  $T_{e,M}(0)\oplus T_{e,M}(\shf)$ is also irreducible as a $T_{e,V}(0)\oplus T_{e,V}(\shf)$-module.
  \\
  (2) The space $T_{e,M}(\sfr{1}{16})$ forms a $\Z_2$-twisted 
  $T_{e,V}(0)\oplus T_{e,V}(\shf)$-module such that 
  $M^{\la -\tau_e\ra}$ as a $V^{\la \tau_e\ra}$-module is isomorphic to a tensor product 
  of a $\Z_2$-twisted $L(\shf,0)\oplus L(\shf,\shf)$-module $L(\shf,\sfr{1}{16})^+$ and 
  the $\Z_2$-twisted $T_{e,V}(0)\oplus T_{e,V}(\shf)$-module $T_{e,M}(\sfr{1}{16})$.
  If $M^{\la -\tau_e\ra}$ is an irreducible $V^{\la \tau_e\ra}$-module, then 
  $T_{e,M}(\sfr{1}{16})$ is also irreducible as a $\Z_2$-twisted 
  $T_{e,V}(0)\oplus T_{e,V}(\shf)$-module.
\end{prop}

\pf
The proof for the existences of structures of modules is similar to that of 
Proposition \ref{prop:3.1} (see Theorem 2.2 of \cite{Y1}).
The irreducibility is clear and follows from the fusion rules of $L(\shf,0)$-modules.
\qed

\begin{rem}\label{rem:3.3}
  If $V$ is simple and both $V_e[\shf]$ and $V_e[\sfr{1}{16}]$ are non-zero, then 
  the fusion rules of $L(\shf,0)$-modules guarantee that 
  $T_{e,M}(h)$ is also non-zero for all non-zero $V$-modules and for $h=0,1/2,1/16$.
  On the other hand, if $V$ is simple, $V_e[\sfr{1}{2}] \ne 0$ and $V_e[\sfr{1}{16}]=0$, 
  then again by the fusion rules of $L(\shf,0)$-modules one of $T_{e,M}(0)\oplus T_{e,M}(\shf)$ 
  or $T_{e,M}(\sfr{1}{16})$ is zero for an irreducible $V$-module $M$. 
\end{rem}

Suppose $e$ is an Ising vector of $V$ such that $T_{e,V}(\shf)\ne 0$.
By Proposition \ref{prop:3.2}, there is a correspondence between 
$V^{\la \tau_e\ra}$-modules and untwisted and $\Z_2$-twisted 
$T_{e,V}(0)\oplus T_{e,V}(\shf)$-modules via \eqref{eq:3.9} and \eqref{eq:3.10}.
It is shown in \cite{Mi3} that $V$ is $C_2$-cofinite if and only if 
$T_{e,V}(0)$ is.
Therefore, we have the following theorem.

\begin{thm}[cf.~\cite{Mi3}]\label{thm:3.4}
  Suppose $e$ is an Ising vector of $V$ such that $T_{e,V}(\shf)\ne 0$.
  \\
  (1)  $V$ is $C_2$-cofinite if and only if $T_{e,V}(0)$ is $C_2$-cofinite.
  \\
  (2)~ $V^{\la \tau_e\ra}$ is rational if and only if 
  $T_{e,V}(0)\oplus T_{e,V}(\shf)$ is rational and $\Z_2$-rational.
\end{thm}

\subsection{A characterization by 3-dimensional Griess algebras}

In this subsection, we will prove the following theorem.

\begin{thm}\label{thm:3.5}
  Let $(V,\w)$ be a VOA of OZ-type.
  Suppose the following.
  \\
  (1) The central charge $c_V$ of $V$ is not equal to $1/2$.
  \\
  (2) $V$ has an Ising vector $e$.
  \\
  (3) There exists a 3-dimensional subalgebra $B=\C \w + \C e + \C x$ of the Griess algebra 
  of $V$ such that $2e_{(1)}x=x$ and $(x\mydiv x)$ is non-zero.
  \\
  Suppose that  $W=\la B\ra$ is a full sub VOA of $V$.
  Then the commutant superalgebra 
  \[
    T_{e,W}(0)\oplus T_{e,W}(\shf)
  \] 
  is isomorphic to an $N=1$ Virasoro SVOA with the conformal vector $\w -e$.
\end{thm}

\pf
Set $f:=\w-e$.
Since $V$ is of OZ-type, $e_{(2)}\w=0$ so that $e$ and $f$ are mutually orthogonal 
Virasoro vectors by Theorem 5.1 of \cite{FZ}.
The central charge of $f$ is $c_f:=2(f\mydiv f)=c_V-1/2\ne 0$.
By a normalization, we may assume that $(x\mydiv x)=2c_f/3$.
Since $x$ is a highest weight vector for $\la e\ra\tensor \la f\ra$ with highest weight 
$(1/2,3/2)$, we have 
\begin{equation}\label{eq:3.14}
  \l( x_{(1)}x\mydiv e\r) = \l( x\mydiv x_{(1)}e\r) =\dfr{(x\mydiv x)}{2}=\dfr{c_f}{3},~~~
  \l( x_{(1)}x\mydiv f\r) = \l( x\mydiv x_{(1)}f\r) =\dfr{3(x\mydiv x)}{2}=c_f.
\end{equation}
By assumption, $\o(e)$ acts on $B$ semisimply with eigenvalues 0, 1/2 and 2.
Thus $e$ is of $\sigma$-type on $W=\la B\ra$ and defines the involution $\sigma_e$ 
by Theorem \ref{thm:2.3}.
Since $\sigma_e$ negates $x$, we have $x_{(1)}x\in B^{\la \sigma_e\ra}=\C e+\C f$ and 
so we can write $x_{(1)}x=\alpha e +\beta f$ with $\alpha$, $\beta\in \C$.
Then 
\begin{equation}\label{eq:3.15}
\begin{array}{l}
  \l( x_{(1)}x\mydiv e\r) = (\alpha e+\beta f\mydiv e) =\alpha (e\mydiv e)=\dfr{\alpha}{4},
  \medskip\\
  \l( x_{(1)}x\mydiv f\r) = (\alpha e+\beta f\mydiv f) =\beta (f\mydiv f)=\dfr{\beta c_f}{2}.
\end{array}
\end{equation}
By \eqref{eq:3.14} and \eqref{eq:3.15}  we have
\begin{equation}\label{eq:3.16}
  x_{(1)}x= \dfr{4c_f}{3} e+ 2f.
\end{equation}
Decompose $W$ as 
\[
  W=L(\shf,0)\tensor T_{e,W}(0)\oplus L(\shf,\shf)\tensor T_{e,W}(\shf).
\]
Then $T:=T_{e,W}(0)\oplus T_{e,W}(\shf)$ is an SVOA with the conformal vector $f$ 
by Proposition \ref{prop:3.1}.
Since $x$ is a highest weight vector for $\la e\ra$ and $\la f\ra$ with highest weight 1/2 
and 3/2, respectively, there exists a highest weight vector $a\in T_{e,W}(\shf)$ for $\la f\ra$ 
with highest weight 3/2 such that $x=\psi_{-1/2}\vac \tensor a$.
Let $Y_V(\cd,z)$ and $Y_T(\cd,z)$ be vertex operator maps of $V$ and $T$, respectively.
Then $Y_V(x,z)=\psi(z)\tensor Y_T(a,z)$ by Proposition \ref{prop:3.1} and we have 
\begin{equation}\label{eq:3.17}
  x_{(n)}x
  = \vac \tensor a_{(n-1)}a + \dsum_{j\geq 0}\psi_{-j-3/2}\psi_{-1/2}\vac \tensor a_{(n+j+1)}a.
\end{equation}
Since $V$ is of OZ-type, so is $T_{e,V}(0)$ and we have\footnote{%
That $a_{(1)}a=0$ also follows from the skew-symmetry.} 
$a_{(1)}a=0$, $a_{(2)}a\in \C\vac$ and $a_{(n)}a=0$ for $n>2$.
Then by \eqref{eq:3.17}, we have 
\begin{equation}\label{eq:3.18}
  x_{(1)}x
  = \vac \tensor a_{(0)}a + \psi_{-3/2}\psi_{-1/2}\vac \tensor a_{(2)}a .
\end{equation}
Since $\psi_{-3/2}\psi_{-1/2}\vac = 2e$, comparing \eqref{eq:3.16} and \eqref{eq:3.18} we obtain 
\begin{equation}\label{eq:3.19}
  a_{(2)}a = \dfr{2c_f}{3} \vac,~~~ 
  a_{(1)}a = 0,~~~
  a_{(0)}a = 2f.
\end{equation}
Now set $L^f(m)=f_{(m+1)}$ and $G^a(r)=a_{(r+1/2)}$ for $m\in \Z$ and $r\in \Z +1/2$. 
By \eqref{eq:3.19} their commutators are as follows.
\[
\begin{array}{ll}
  {}\l[ L^f(m),G^a(r)\r] 
  &= \l[ f_{(m+1)},a_{(r+1/2)}\r] 
  = \dsum_{i=0}^\infty \binom{m+1}{i} \l( f_{(i)}a\r)_{(m+r+3/2-i)}
  \medskip\\
  &= \l( f_{(0)}a\r)_{(m+r+3/2)}+(m+1)\l( f_{(1)}a\r)_{(m+r+1/2)}
  \medskip\\
  &= -\l( m+r+\dfr{3}{2}\r) a_{(m+r+1/2)}+\dfr{3}{2}(m+1)a_{(m+r+1/2)}
  \medskip\\
  &= \l( \dfr{1}{2}m -r\r) a_{(m+r+1/2)}
  = \l( \dfr{1}{2}m -r\r) G^a(m+r),
\end{array}
\]

\[
\begin{array}{ll}

  \l[ G^a(r),G^a(s)\r]
  &= \l[ a_{(r+1/2)},a_{(s+1/2)}\r]
  = \dsum_{i=0}^\infty \binom{r+1/2}{i} \l( a_{(i)}a\r)_{(r+s+1-i)}
  \medskip\\
  &= \ds \l( a_{(0)}a\r)_{(r+s+1)}
     +\binom{r+1/2}{2} \l( a_{(2)}a\r)_{(r+s-1)}
  \medskip\\
  &= \ds 2f_{(r+s+1)}+\binom{r+1/2}{2}\cd \dfr{2c_f}{3}\cd \vac_{(r+s-1)}
  \medskip\\
  &= 2L^f(r+s)+\delta_{r+s,0}\l( r^2-\dfr{1}{4}\r)\dfr{c_f}{3}.
\end{array}
\]
Therefore, $Y_T(f,z)$ and $Y_T(a,z)$ generate a representation of the 
Neveu-Schwarz algebra $\ns$ on $T$.
Since $W$ is generated by $e$, $f$ and $x=\psi_{-1/2}\vac \tensor a$, 
it follows that $T=T_{e,W}(0)\oplus T_{e,W}(\shf)$ is generated by $f$ and $a$.
Therefore, $T$ is isomorphic to an $N=1$ super Virasoro VOA.
\qed

\section{Extension of a pair of unitary Virasoro VOAs}

Consider the even part 
\begin{equation}\label{eq:4.1}
  A(\shf,c_n^1)
  := L(\shf,0)\tensor L_\ns(c_n^1,0)^{[0]} \oplus L(\shf,\shf)\tensor L_\ns(c_n^1,0)^{[1]} 
\end{equation}
of the tensor product of SVOAs $L(\shf,0)\oplus L(\shf,\shf)$ and 
$L_\ns(c_n^1,0)$ where $c_n^1$ is defined as in \eqref{eq:2.11}.
This VOA is also considered in \cite{A2,LS}.
The VOA $A(\shf,c_n^1)$ inherits the invariant bilinear forms of 
$L(\shf,0)\oplus L(\shf,\shf)$ and $L_\ns(c_n^1,0)$ and has a compact real form.

\subsection{Griess algebra}

Let $e=\fr{1}{2}\psi_{-3/2}\psi_{-1/2}\vac$ and $f=\fr{1}{2}G(-\shf)G(-\sfr{3}{2})\vac$ 
be the conformal vectors of $L(\shf,0)$ and $L_\ns(c_m^1,0)^{[0]}$, respectively, and 
let $x=\sqrt{(n+2)(n+4)}\, \psi_{-1/2}\vac \tensor G(-\sfr{3}{2})\vac$ be the highest weight vector 
of $L(\shf,\shf)\tensor L_\ns(c_n^1,0)^{[1]}$.
Then $A(\shf,c_n^1)$ is of OZ-type and its Griess algebra is 3-dimensional with an orthogonal 
basis $e$, $f$ and $x$ such that 
\begin{equation}\label{eq:4.2}
\begin{array}{l}
  e_{(1)}e=2e,~~~
  e_{(1)}f=0,~~~
  e_{(1)}x=\dfr{1}{2}x,~~~
  f_{(1)}f=2f,~~~
  f_{(1)}x=\dfr{3}{2}x,~~~
  (e\mydiv e)=\dfr{1}{4},
  \medskip\\
  x_{(1)}x=2n(n+6)e+2(n+2)(n+4)f,~~~
  (f\mydiv f)=\dfr{c_n^1}{2},~~~(x\mydiv x)=n(n+6).
\end{array}
\end{equation}
By a direct calculation,  we can classify the Virasoro vectors in $A(\shf,c_n^1)$.

\begin{prop}[\cite{A2,LS}]\label{prop:4.1}
  Let 
  \begin{equation}\label{eq:4.3}
    u := \dfr{1}{2(n+3)}\l(  n e +(n+4) f+x\r) ,\quad \text{ and } \quad  
    v := e+f-u.
  \end{equation}
  (1) $u$ and $v$ are mutually orthogonal Virasoro vectors with central charges 
  $c_n^0$ and $c_{n+1}^0$.
  \\
  (2)The set of Virasoro vectors of $A(\shf,c_n^1)$ is given by 
  $\{ \w, e,f, u,v,\sigma_e u, \sigma_e v\}$.
  \\
  (3) $A(\shf,c_n^1)$ is generated by its Griess algebra.
  \\
  (4)  
  $\aut(A(\shf,c_n^1))=\la \sigma_e\ra$ if $n>1$ 
  and $\aut(A(\shf,c_1^1))=\la \sigma_e,\sigma_{u}\ra\cong \mathrm{S}_3$.
\end{prop}

\pf
(1): It is straightforward to verify that
$u$ and $v$ are mutually orthogonal Virasoro vectors and $\w=u+v$ is a Virasoro frame of 
$A(\shf,c_n^1)$ (cf.~\cite{A2}).
\\
(2): The solutions of the quadratic equation $y^2=2y$ in the Griess algebra provide 
a complete list of Virasoro vectors in $A(\shf,c_n^1)$ which is as in the assertion.
\\
(3): Let $V$ be a subalgebra of $A(\shf,c_n^1)$ generated by the Griess algebra.
Then $V$ satisfies the conditions in Theorem \ref{thm:3.5} and it follows that 
$V=A(\shf,c_n^1)$.
\\
(4): If $n>1$ then $e$ is the unique Ising vector of $A(\shf,c_n^1)$ and 
$\sigma_e$ is the unique non-trivial automorphism of the Griess algebra.
Since $A(\shf,c_n^1)$ is generated by its Griess algebra, we have 
$\aut(A(\shf,c_n^1))=\la \sigma_e\ra$.
If $n=1$, then $A(\shf,c_n^1)$ has three Ising vectors $e$, $u$ and $\sigma_e u$ of 
$\sigma$-types and $\sigma$-involutions associated to these Ising vectors generate 
$\mathrm{S}_3$.
\qed
\medskip

Since $A(\shf,c_n^1)$ has a compact real form, it contains a full sub VOA 
$\la u\ra\tensor \la v\ra\cong L(c_n^0,0)\tensor L(c_{n+1}^0,0)$.
By \eqref{eq:4.2} and \eqref{eq:4.3}, the Griess algebra is spanned by 
$e$, $u$ and $v$ where their multiplications are as follows.
\begin{equation}\label{eq:4.4}
\begin{array}{l}
  e_{(1)}u=\dfr{n+1}{n+3}e+\dfr{n+2}{4(n+3)}u-\dfr{n+4}{4(n+3)}v,
  \medskip\\
  e_{(1)} v=\dfr{n+5}{n+3}e-\dfr{n+2}{4(n+3)}u+\dfr{n+4}{4(n+3)}v.
\end{array}
\end{equation}
Clearly, the relations above uniquely determines the Griess algebra 
so that $A(\shf,c_n^1)$ is also characterized by the structure in \eqref{eq:4.4} 
thanks to Theorem \ref{thm:3.5}.

\begin{prop}\label{prop:4.2}
Let $V_\R$ be a compact VOA of OZ-type.
Suppose $e$, $u$ and $v$ are simple $c=c_1^0$, $c=c_n^0$ and $c=c_{n+1}^0$ Virasoro vectors 
of $V_\R$, respectively, such that $u$ and $v$ are mutually orthogonal and satisfy \eqref{eq:4.4}.
Then the sub VOA generated by $e$, $u$, and $v$ is isomorphic to $A(\shf,c_n^1)$.
\end{prop}

\pf
Suppose $u_{(1)}v=0$ and \eqref{eq:4.4}.
The inner products $(e\mydiv u)$ and $(e\mydiv v)$ are uniquely determined by 
the invariance property  
$(e_{(1)}u|u)=(e|u_{(1)}u)=2(e|u)$ and $(e_{(1)}v|v)=(e|v_{(1)}v)=2(e|v)$. 
Then by change of basis, we recover the relations \eqref{eq:4.2} and hence the subalgebra 
generated by $u$, $v$ and $e$ is isomorphic to $A(\shf,c_n^1)$ by Theorem \ref{thm:3.5}.
\qed
\medskip

Set 
\begin{equation}\label{eq:4.5}
  w = 3n(n+6)e-(n+2)(n+4)f+3x.
\end{equation}
Then $w$ is a highest weight vector for $\la u\ra\tensor \la v\ra$ 
with the highest weight 
\[
  (h_{1,3}^{(n)},h_{3,1}^{(n+1)})=\l(\dfr{n+1}{n+3},\dfr{n+5}{n+3}\r) .
\]
Therefore, $A(\shf,c_n^1)$ contains $L(c_n^0,h_{1,3}^{(n)})\tensor L(c_{n+1}^0,h_{3,1}^{(n+1)})$ 
as a $\la u\ra\tensor \la v\ra$-submodule.
A complete decomposition will be given in the next subsection.

\begin{thm}[cf.~\cite{A2}]\label{thm:4.3}~
  \\
  (1) $A(\shf,c_n^1)$ is rational and $C_2$-cofinite.
  \\
  (2) The even part $L_\ns(c_n^1,0)^{[0]}$ is $C_2$-cofinite.
\end{thm}

\pf
The rationality of $A(\shf,c_n^1)$ follows from Theorem \ref{thm:2.4} and (2) of 
Theorem \ref{thm:3.4}.
Since $A(\shf,c_n^1)$ has a $C_2$-cofinite full sub VOA $L(c_n^0,0)\tensor L(c_{n+1}^0,0)$, 
it is also $C_2$-cofinite (cf.~\cite{ABD,A2}), 
and the $C_2$-cofiniteness of $L_\ns(c_n^1,0)^{[0]}$ follows from (1) of Theorem \ref{thm:3.4}.
\qed

\subsection{Modules}

Let $n$ be a positive integer and let 
$L_{\hat{\mathfrak{sl}}_2}(n,j)=L_{\hat{\mathfrak{sl}}_2}((n-j)\Lambda_0+j\Lambda_1)$ be the 
level $n$ integrable highest weight $\hat{\mathfrak{sl}}_2$-module with highest weight 
$(n-j)\Lambda_0+j\Lambda_1$, $0\leq j\leq n$.
By \cite{GKO}, $L_{\hat{\mathfrak{sl}}_2}(1,0)\tensor L_{\hat{\mathfrak{sl}}_2}(n,0)$ contains a 
full sub VOA $L(c_n^0,0)\tensor L_{\hat{\mathfrak{sl}}_2}(n+1,0)$ and we have the following 
decompositions for $i=0,1$ and $0\leq j\leq n$.
\begin{equation}\label{eq:4.6}
  L_{\hat{\mathfrak{sl}}_2}(1,i)\tensor L_{\hat{\mathfrak{sl}}_2}(n,j)
  =\bigoplus_{0\leq k\leq n+1 \atop k\equiv i+j (2)} 
  L(c_n^0,h_{j+1,k+1}^{(n)})\tensor L_{\hat{\mathfrak{sl}}_2}(n+1,k).
\end{equation}
By \cite{Li2}, the affine VOA $L_{\hat{\mathfrak{sl}}_2}(2,0)$ admits an extension to a 
simple SVOA $L_{\hat{\mathfrak{sl}}_2}(2,0)\oplus L_{\hat{\mathfrak{sl}}_2}(2,2)$ 
by a simple current module $L_{\hat{\mathfrak{sl}}_2}(2,2)$.
The classifications of irreducible untwisted and $\Z_2$-twisted 
$L_{\hat{\mathfrak{sl}}_2}(2,0)\oplus L_{\hat{\mathfrak{sl}}_2}(2,2)$-modules are established
in (loc.~cit.).
The adjoint module is the unique irreducible untwisted 
$L_{\hat{\mathfrak{sl}}_2}(2,0)\oplus L_{\hat{\mathfrak{sl}}_2}(2,2)$-module and 
there exist two inequivalent structures $L_{\hat{\mathfrak{sl}}_2}(2,1)^\pm$ of irreducible 
$\Z_2$-twisted $L_{\hat{\mathfrak{sl}}_2}(2,0)\oplus L_{\hat{\mathfrak{sl}}_2}(2,2)$-modules on 
$L_{\hat{\mathfrak{sl}}_2}(2,1)$ which are mutually $\Z_2$-conjugate to each other. 
It is shown in \cite{GKO} that a tensor product
\[
  (L_{\hat{\mathfrak{sl}}_2}(2,0)\oplus L_{\hat{\mathfrak{sl}}_2}(2,2))\tensor 
  L_{\hat{\mathfrak{sl}}_2}(n,0)
\]
contains a full sub SVOA $L_\ns(c_n^1,0)\tensor L_{\hat{\mathfrak{sl}}_2}(n+2,0)$ and 
we have the following decompositions for $i=0$, $1$ and $0\leq j\leq n$.
\begin{equation}\label{eq:4.7}
\begin{array}{l}
  \ds L_{\hat{\mathfrak{sl}}_2}(2,2i)\tensor L_{\hat{\mathfrak{sl}}_2}(n,j)
  =\bigoplus_{0\leq k\leq n+2 \atop k\equiv j\,(2)} 
   L_{\mathrm{NS}}(c_n^1,h_{j+1,k+1,0}^{(n)})^{[i+\fr{j-k}{2}]}
   \tensor L_{\hat{\mathfrak{sl}}_2}(n+2,k), 
  \medskip\\
  \l( L_{\hat{\mathfrak{sl}}_2}(2,1)^+\oplus L_{\hat{\mathfrak{sl}}_2}(2,1)^-\r) 
  \tensor L_{\hat{\mathfrak{sl}}_2}(n,j)
  \medskip\\
  \ds = \bigoplus_{0\leq k\leq n+2 \atop k\equiv j+1\,(2)} 
    \l( L_{\ram}(c_n^1,\Delta_{j+1,k+1}^{(n)}) \oplus 
    L_{\ram}(c_n^1,-\Delta_{j+1,k+1}^{(n)})\r) 
   \tensor L_{\hat{\mathfrak{sl}}_2}(n+2,k) .
\end{array}
\end{equation}
As we have seen, 
$A(\shf,c_n^1)=L(\shf,0)\tensor L_{\mathrm{NS}}(c_n^1,0)^{[0]}\oplus L(\shf,\shf)\tensor 
L_{\mathrm{NS}}(c_n^1,0)^{[1]}$ contains a full sub VOA $L(c_n^0,0)\tensor L(c_{n+1}^0,0)$.
We consider the decompositions of irreducible $A(\shf,c_n^1)$-modules 
as $L(c_n^0,0)\tensor L(c_{n+1}^0,0)$-modules.

First, we label the irreducible $A(\shf,c_n^1)$-modules as follows.
\begin{equation}\label{eq:4.8}
\begin{array}{l}
  M(0,h_{r,s,0}^{(n)})
  = L(\shf,0)\tensor L_\ns(c_n^1,h_{r,s,0}^{(n)})^{[0]}
  \oplus L(\shf,\shf)\tensor L_\ns(c_n^1,h_{r,s,0}^{(n)})^{[1]},
  \medskip\\
  M(\shf,h_{r,s,0}^{(n)})
  = L(\shf,\shf)\tensor L_\ns(c_n^1,h_{r,s,0}^{(n)})^{[0]}
  \oplus L(\shf,0)\tensor L_\ns(c_n^1,h_{r,s,0}^{(n)})^{[1]},
  \medskip\\
  M(\sfr{1}{16},\Delta_{r,s}^{(n)})
  =L(\shf,\sfr{1}{16})^+ \tensor L_{\ram}(c_n^1,\Delta_{r,s}^{(n)}) ,
\end{array}
\end{equation}
where $h_{r,s,p}^{(n)}$ and $\Delta_{r,s}^{(n)}$ with $1\leq r\leq n+1$, $1\leq s\leq n+3$,
$p=0$, $1$ and $r-s\equiv 2p \mod 2$, are as in \eqref{eq:2.11}.
The zero-mode $\o(\psi_{-1/2}\vac \tensor G(-\sfr{3}{2})\vac)=\phi_0\tensor G(0)$ acts on 
the top level of $M(\sfr{1}{16},\Delta_{r,s}^{(n)})$ 
by 
\begin{equation}\label{eq:4.9}
  \dfr{1}{\sqrt{2}}\cd \Delta_{r,s}^{(n)}
  =\dfr{s(n+2)-r(n+4)}{4\sqrt{(n+2)(n+4)}} . 
\end{equation}
Note that $\Delta_{n+2-r,n+4-s}^{(n)}=-\Delta_{r,s}^{(n)}$ so that 
\[ 
  L(\sfr{1}{16},0)^- \tensor L_{\ram}(c_n^1,\Delta_{r,s}^{(n)})
  \cong 
  L(\sfr{1}{16},0)^+ \tensor L_{\ram}(c_n^1,\Delta_{n+2-r,n+4-s}^{(n)})
\]
as $A(\shf,c_n^1)$-modules and $M(\sfr{1}{16},\pm \Delta_{r,s}^{(n)})$ are mutually 
$\sigma_e$-conjugate while $M(\shf,h_{r,s,0}^{(n)})$ are $\sigma_e$-invariant.
The next theorem follows from Proposition \ref{prop:3.2}.

\begin{thm}\label{thm:4.4}
  The set of irreducible $A(\shf,c_n^1)$-modules is given by the list \eqref{eq:4.8}.
\end{thm}

We decompose an irreducible $A(\shf,c_n^1)$-module into a direct sum of irreducible 
modules over $L(c_n^0,0)\tensor L(c_{n+1}^0,0)$.
By \eqref{eq:4.6}, we have
\begin{equation}\label{eq:4.10}
\begin{array}{l}
  L_{\hat{\mathfrak{sl}}_2}(1,0)\tensor L_{\hat{\mathfrak{sl}}_2}(1,0)  
  = L(\shf,0)\tensor L_{\hat{\mathfrak{sl}}_2}(2,0) \oplus 
  L(\shf,\shf)\tensor L_{\hat{\mathfrak{sl}}_2}(2,2), 
  \medskip\\
  L_{\hat{\mathfrak{sl}}_2}(1,1)\tensor L_{\hat{\mathfrak{sl}}_2}(1,1)  
  = L(\shf,0)\tensor L_{\hat{\mathfrak{sl}}_2}(2,2) \oplus 
  L(\shf,\shf)\tensor L_{\hat{\mathfrak{sl}}_2}(2,0), 
  \medskip\\
  L_{\hat{\mathfrak{sl}}_2}(1,0) \tensor L_{\hat{\mathfrak{sl}}_2}(1,1)  
  \oplus L_{\hat{\mathfrak{sl}}_2}(1,1)\tensor L_{\hat{\mathfrak{sl}}_2}(1,0)  
  \medskip\\
  = \l( L(\shf,\sfr{1}{16}) \tensor L_{\hat{\mathfrak{sl}}_2}(2,1)\r)^+ 
  \oplus \l( L(\shf,\sfr{1}{16}\tensor L_{\hat{\mathfrak{sl}}_2}(2,1)\r)^-.
\end{array}
\end{equation}
Plugging \eqref{eq:4.10} into \eqref{eq:4.7}, we obtain the following.

\begin{prop}[\cite{A2,LS}]\label{prop:4.5}
As $L(c_n^0,0)\tensor L(c_{n+1}^0,0)$-modules, we have the following decompositions.
\[
\begin{array}{l}
  M(\varepsilon,h_{r,s,0}^{(n)})
  = \ds\bigoplus_{1\leq j\leq n+2 \atop j\equiv (r+s)/2+2\varepsilon \,(2)} 
    L(c_n^0,h_{r,j}^{(n)}) \tensor L(c_{n+1}^0,h_{j,s}^{(n+1)}),~~~
    \varepsilon =0,\dfr{1}{2},
  \medskip\\
  M(\sfr{1}{16},\pm \Delta_{r,s}^{(n)})
  = \ds\bigoplus_{1\leq j\leq n+2 \atop j\equiv (r+s \pm 1)/2\, (2)}
    L(c_n^0,h_{r,j}^{(n)}) \tensor L(c_{n+1}^0,h_{j,s}^{(n+1)}).
\end{array}
\]
\end{prop}

\pf
The decomposition of $M(\varepsilon,h_{r,s,0}^{(n)})$ with $\varepsilon = 0,1/2$ 
is straightforward (cf.~\cite{A2}).
By \eqref{eq:4.6} and \eqref{eq:4.7} we obtain
\[
  M(\sfr{1}{16},\Delta_{r,s}^{(n)}) \oplus M(\sfr{1}{16},-\Delta_{r,s}^{(n)})
  = \bigoplus_{1\leq j\leq n+2} L(c_n^0,h_{r,j}^{(n)})\tensor L(c_{n+1}^0,h_{j,s}^{(n+1)}).
\]
By \eqref{eq:4.3} and \eqref{eq:4.9}, the top levels of 
$M(\sfr{1}{16},\pm \Delta_{r,s}^{(n)})$ contain those of 
\[ 
  L(c_n^0,h_{r,(r+s\pm 1)/2}^{(n)})\tensor L(c_{n+1}^0,h_{(r+s\pm 1)/2,s}^{(n+1)}), 
\]
respectively. 
Since 
\[
  A(\shf,c_n^1) = M(0,h_{1,1,0}^{(n)})
  =\bigoplus_{1\leq j\leq n+2 \atop j\equiv 1\, (2)} L(c_n^0,h_{1,j}^{(n)})\tensor L(c_n^0,h_{j,1}^{(n+1)}),
\]
the decompositions of $M(\sfr{1}{16},\pm \Delta_{r,s}^{(n)})$ are determined by the fusion rules 
of $L(c_n^0,0)$ and $L(c_{n+1}^0,0)$-modules in \eqref{eq:2.3} as in the assertion.
\qed

\subsection{Automorphisms}

We will determine the group generated by Miyamoto involutions of $c=c_n^0$ and $c=c_{n+1}^0$ 
Virasoro vectors in  $A(\shf,c_n^1)$.

\begin{thm}\label{thm:4.6}
  Suppose a VOA $V$ contains a sub VOA $U$ isomorphic to $A(\shf,c_n^1)$.
  Let $e$, $u$ and $v$ be $c=c_1^0$, $c=c_n^0$ and $c=c_{n+1}^0$ 
  Virasoro vectors of $U$ given  by \eqref{eq:4.3}, respectively. 
  Then the following hold.
  \\
  (1) $[\tau_u,\tau_v]=[\tau_{\sigma_eu},\tau_{\sigma_e v}]=1$ in $\aut(V)$ and 
  $\tau_e$ centralizes $\la \tau_u,\tau_v,\tau_{\sigma_e u},\tau_{\sigma_e v}\ra$.
  \\
  (2) If $n$ is even then $\tau_u\tau_v=\tau_{\sigma_eu}\tau_{\sigma_ev}=\tau_e$, 
  $\tau_u=\tau_{\sigma_e u}$ and $\tau_v=\tau_{\sigma_e v}$ in $\aut(V)$.
  \\
  (3) If $n$ is odd then $\tau_u=\tau_v$, $\tau_{\sigma_eu}=\tau_{\sigma_e v}$ and 
  $\tau_u\tau_{\sigma_e u}=\tau_v\tau_{\sigma_e v}=\tau_e$ in $\aut(V)$.
  \\
  (4) $\la \tau_e,\tau_u,\tau_v,\tau_{\sigma_eu},\tau_{\sigma_ev}\ra$ is 
  an elementary abelian 2-group of rank at most 2.
\end{thm}

\pf
Since $\w=u+v=\sigma_eu+\sigma_ev$ are Virasoro frames of $U$, 
we have $[\tau_u,\tau_v]=[\tau_{\sigma_eu},\tau_{\sigma_e v}]=1$ in $\aut(V)$.
Since $e$ is of $\sigma$-type on $U$, $\tau_e$ is trivial on $U$.
Then it follows from $\tau_y=\tau_{\tau_e y}=\tau_e\tau_y\tau_e$ for 
$y\in \{ u,v,\sigma_eu,\sigma_ev\}$ that $\tau_e$ centralizes 
$\la \tau_u,\tau_v,\tau_{\sigma_e u},\tau_{\sigma_e v}\ra$. 
By (1) of Theorem \ref{thm:4.3}, $V$ is a direct sum of irreducible $A(\shf,c_n^1)$-submodules.
By definition, Miyamoto involutions preserve each irreducible $A(\shf,c_n^1)$-modules.
If $n$ is even, it follows from the decompositions in Proposition \ref{prop:4.5} that 
$\tau_u$ and $\tau_v$ act on $M(\varepsilon,h_{r,s,0}^{(n)})$, 
$\varepsilon=0$, $1/2$ and $M(\sfr{1}{16},\Delta_{r,s}^{(n)})$ as $(-1)^{r+1}$ and $(-1)^{s+1}$, 
respectively.
Then the product $\tau_u\tau_v=(-1)^{r+s}$ is trivial on 
$M(\varepsilon,h_{r,s,0}^{(n)})$ and is equal to $-1$ on 
$M(\sfr{1}{16},\Delta_{r,s}^{(n)})$ since $r\equiv s \mod 2$ for the NS-sectors and 
$r\not\equiv s\mod 2$ for the R-sectors.
On the other hand, $\tau_e$ is trivial on $M(\varepsilon,h_{r,s,0}^{(n)})$, $\varepsilon=0$, $1/2$ 
and acts as $-1$ on $M(\sfr{1}{16},\Delta_{r,s}^{(n)})$.
Thus $\tau_u\tau_v=\tau_e$ in $\aut(V)$.
Since $M(\varepsilon,h_{r,s,0}^{(n)})$ is $\sigma_e$-invariant and 
$M(\sfr{1}{16},\pm \Delta_{r,s}^{(n)})$ are mutually $\sigma_e$-conjugate, 
we have $\tau_u=\tau_{\sigma_eu}$ and $\tau_v=\tau_{\sigma_e v}$ in $\aut(V)$.

If $n$ is odd then both $\tau_u$ and $\tau_v$ act on 
$L(c_n^0,h_{r,j}^{(n)})\tensor L(c_{n+1}^0,h_{j,s}^{(n+1)})$ by $(-1)^{j+1}$ and  
it follows from the decompositions in Proposition \ref{prop:4.5} that 
$\tau_u=\tau_v$ in $\aut(V)$.
Each summand $L(c_n^0,h_{r,j}^{(n)})\tensor L(c_{n+1}^0,h_{j,s}^{(n+1)})$ of 
$M(\varepsilon,h_{r,s,0}^{(n)})$ with $\varepsilon=0$, $1/2$ satisfies 
$j\equiv (r+s)/2+2\varepsilon \mod 2$ and hence $\tau_u=\tau_v$ acts as $(-1)^{(r+s)/2+2\varepsilon}$ 
on $M(\varepsilon,h_{r,s,0}^{(n)})$.
Since $M(\varepsilon,h_{r,s,0}^{(n)})$ is $\sigma_e$-stable, its decomposition 
with respect to $\la \sigma_e u\ra \tensor \la \sigma_e v\ra$ is isomorphic to 
that with respect to $\la u\ra \tensor \la v\ra$.
Therefore $\tau_{\sigma_e u}=\tau_{\sigma_e v}$ also satisfies 
$\tau_{\sigma_e u}=\tau_{\sigma_e v}=(-1)^{(r+s)/2+2\varepsilon}$ on $M(\varepsilon,h_{r,s,0}^{(n)})$.
Thus $\tau_u\tau_{\sigma_e u}=\tau_v\tau_{\sigma_e v}=1=\tau_e$ on $M(\varepsilon,h_{r,s,0}^{(n)})$.
On the other hand, each summand $L(c_n^0,h_{r,j}^{(n)})\tensor L(c_{n+1}^0,h_{j,s}^{(n+1)})$ of 
$M(\sfr{1}{16},\Delta_{r,s}^{(n)})$ satisfies $j\equiv (r+s+1)/2 \mod 2$ and hence 
$\tau_u=\tau_v$ acts as $(-1)^{(r+s+1)/2}$ on $M(\sfr{1}{16},\Delta_{r,s}^{(n)})$.
Since the decomposition of $M(\sfr{1}{16},\Delta_{r,s}^{(n)})$ as a 
$\la \sigma_e u\ra\tensor \la \sigma_e v\ra$-module is isomorphic to 
that of $M(\sfr{1}{16},-\Delta_{r,s}^{(n)})$ as a $\la u\ra\tensor \la v\ra$-module, 
we have $\tau_{\sigma_e u}=\tau_{\sigma_e v}=(-1)^{(r+s-1)/2}$ on $M(\sfr{1}{16},\Delta_{r,s}^{(n)})$.
Therefore $\tau_u\tau_{\sigma_e u}=\tau_v\tau_{\sigma_e v}=(-1)^{r+s}=-1=\tau_e$ on 
$M(\sfr{1}{16},\Delta_{r,s}^{(n)})$.
This completes the proof.
\qed

\begin{rem}\label{rem:4.7}
  (3) of Theorem \ref{thm:4.6} is also proved in \cite{HLY1} in the case of $n=1$. 
\end{rem}

\begin{thm}\label{thm:4.8}
  Suppose $n$ is odd and a VOA $V$ contains a sub VOA $U$ isomorphic to $A(\shf,c_n^1)$.
  Let $e$, $u$ and $v$ be $c=c_1^0$, $c=c_n^0$ and $c=c_{n+1}^0$ Virasoro vectors of $U$ given 
  by \eqref{eq:4.3}, respectively. 
  Then $v$ is of $\sigma$-type on the commutant $\com_V \la u\ra$.
  Moreover, $\sigma_v$ and $\tau_e$ define the same automorphism of $\com_V \la u\ra$.
\end{thm}

\pf
Let $X$ be an irreducible $\la v\ra$-submodule of $V$.
Then $X\subset \com_V\la u\ra$ if and only if there exists an irreducible 
$A(\shf,c_n^1)$-submodule $M$ of $V$ containing a $\la u\ra\tensor \la v\ra$-submodule 
isomorphic to $L(c_n^0,0)\tensor X$.
By Theorem \ref{thm:4.4} and Proposition \ref{prop:4.5},  
$X\cong L(c_{n+1}^0,h_{1,s}^{(n+1)})$ and 
$M\cong M(\varepsilon,h_{1,s,0}^{(n)})$ with $s\equiv 4\varepsilon+1 \mod 4$ for 
$\varepsilon =0$, $1/2$ or 
$M\cong M(\sfr{1}{16},\pm \Delta_{1,s}^{(n)})$ with $s\equiv 3\pm 1 \mod 4$.
Therefore, $v$ is of $\sigma$-type on $\com_V \la u\ra$.
From the above possibility of $X$ and $M$, both $\sigma_v$ and $\tau_e$ define 
the same automorphism on the commutant $\com_V\la u\ra$.
This completes the proof.
\qed

\begin{rem}\label{rem:4.9}
  Theorem \ref{thm:4.8} is also proved in \cite{HLY1} in the case of $n=1$ 
  and in \cite{HLY2} in the case of $n=3$.
\end{rem}

\appendix
\def\thesection{\Alph{section}}
\section{Appendix}

In this appendix we prove the $\Z_2$-rationality of $L_\ns(c_n^1,0)$.
The classification of irreducible $\Z_2$-twisted $L_\ns(c_n^1,0)$-modules is accomplished 
in \cite{M} in the category of superspaces.
Our argument is almost the same as in (loc.\ cit.) but we do not assume the superspace 
structure on $\Z_2$-twisted modules.
First we recall the $\Z_2$-twisted Zhu algebra of an SVOA.
Let $V=V^{[0]}\oplus V^{[1]}$ be an SVOA such that $V^{[i]}$ has $(\Z+i/2)$-grading.
For homogeneous $a$ and $b\in V$, we define
\begin{equation}\label{eq:a.1}
\begin{array}{l}
  a{\ds \mystar_\tw} b
  := \res_z Y(a,z)b\dfr{(1+z)^{\wt (a)}}{z}
  = \dsum_{i=0}^\infty \binom{\wt(a)}{i}a_{(i-1)}b,
  \medskip\\
  a{\ds \mycirc_\tw} b:= \res_z Y(a,z)b\dfr{(1+z)^{\wt (a)}}{z^2}
  = \dsum_{i=0}^\infty \binom{\wt(a)}{i}a_{(i-2)}b,
\end{array}
\end{equation}
and extend bilinearly.
We then set
\begin{equation}\label{eq:a.2}
  A_\tw(V):=V/O_\tw(V),~~~~~
  O_\tw(V):=\Span_\C \{ a{\ds \mycirc_\tw} b \mid a,b\in V\} .
\end{equation}
We denote the class $a+O_\tw(V)$ of $a\in V$ in $A_\tw(V)$ by $[a]$.
It is shown in \cite{DZ} that $A_\tw(V)$ equipped with the product $\ds \mystar_\tw$ 
in \eqref{eq:a.1} forms a unital associative algebra such that $[\vac]$ is the unit and 
$[\w]$ is in the center.
The twisted Zhu algebra $A_\tw(V)$ determines irreducible $\Z_2$-twisted representations.

\begin{thm}[\cite{Z,DZ}]\label{thm:a.1}~
  \\
  Let $V=V^{[0]}\oplus V^{[1]}$ be an SVOA such that $V^{[i]}$ has $(\Z+i/2)$-grading.
  \\
  (1) Let $M$ be a $\Z_2$-twisted $V$-module and $\Omega(M)$ its top level. 
  Then the zero-mode $\o(a)=a_{(\wt(a)-1)}$ defines a representation of $A_\tw(V)$ on $\Omega(M)$.
  \\
  (2) Let $N$ be an irreducible $A_\tw(V)$-module.
  Then there exists the unique irreducible $\Z_2$-twisted $V$-module $\tilde{N}$ such that 
  its top level is isomorphic to $N$ as $A_\tw(V)$-modules.
  \\
  (3) There is a one-to-one correspondence between irreducible $\Z_2$-twisted $V$-modules and 
  irreducible $A_\tw(V)$-modules.
\end{thm}

We classify irreducible $\Z_2$-twisted modules over the $N=1$ Virasoro SVOA 
$L_\mathrm{NS}(c_n^1,0)$ based on Theorem \ref{thm:a.1}.
We first consider the $\Z_2$-twisted Zhu algebra $A_\tw\!\l(\ol{M}_\mathrm{NS}(c,0)\r)$ 
of the universal $N=1$ Virasoro SVOA 
$\ol{M}_\mathrm{NS}(c,0)=M_\mathrm{NS}(c,0)/\la G(-\sfr{1}{2})\vac)\ra$.
Note that $\ol{M}_\mathrm{NS}(c,0)$ has a linear basis 
\begin{equation}\label{eq:a.3}
  L(-n_1)\cds L(-n_i)G(-r_1)\cds G(-r_j)\vac,~~
  n_1\geq \cds \geq n_i\geq 2,~~
  r_1>\cds >r_j\geq \dfr{3}{2}.
\end{equation}
Images of Virasoro descendants in $A_\tw\!\l(\ol{M}_\mathrm{NS}(c,0)\r)$ are easy to compute.

\begin{lem}\label{lem:a.2}
  $[L(-n)a]=(-1)^n (n-1)[a] \,{\ds \mystar_\tw}\, [\w]+(-1)^n [L(0)a]$ for $n\geq 1$. 
\end{lem}

\pf 
See (4.2) of \cite{W}.
\qed
\medskip

For odd elements we have the following recursion.

\begin{lem}\label{lem:a.3}
  $[G(-r)a]=-\dsum_{s<r}\binom{3/2}{r-s}[G(-s)a]$ for $r\geq 5/2$.
\end{lem}

\pf
By Lemma 2.1.2 of \cite{Z}, for any $n\geq 0$ one has 
\[
  \res_zY(\tau,z)a\dfr{(1+z)^{3/2}}{z^{2+n}}
  = \dsum_{i\geq 0}\binom{3/2}{i} \tau_{(-2-n+i)}a 
  \in O_\tw\l(\ol{M}_{\mathrm{NS}}(c,0)\r) .
\]
Noting $\tau_{(-2-n+i)}=G(-\sfr{5}{2}-n+i)$, we obtain the lemma.
\qed

\begin{lem}\label{lem:a.4}
  $[\tau]^2=[\w]-\dfr{c}{24}[\vac]$.
\end{lem}

\pf
By definition, one has
\[
\begin{array}{ll}
  \tau{\ds \mystar_\tw}\tau
  &=\dsum_{i\geq 0}\binom{3/2}{i}\tau_{(i-1)}\tau
  =\dsum_{i\geq 0}\binom{3/2}{i}G(i-\sfr{3}{2})G(-\sfr{3}{2})\vac
  \medskip\\
  &= \ds G(-\sfr{3}{2})^2\vac +\binom{3/2}{i} G(-\sfr{1}{2})G(-\sfr{3}{2})\vac 
  +\binom{3/2}{3}G(\sfr{3}{2})G(-\sfr{3}{2})\vac
  \medskip\\
  &=L(-3)\vac +3L(-2)\vac -\dfr{c}{24}\vac .
\end{array}
\]
Then by Lemma \ref{lem:a.2}, one obtains
\[
  [\tau]^2
  = [\tau{\ds \mystar_\tw} \tau]
  = [L(-3)\vac] +3[L(-2)\vac] -\dfr{c}{24}[\vac]
  = [\w]-\dfr{c}{24}[\vac].
\]
This completes the proof.
\qed

\begin{prop}\label{prop:a.5}
  The $\Z_2$-twisted Zhu algebra $A_\tw\!\l(\ol{M}_\mathrm{NS}(c,0)\r)$ is isomorphic to 
  the polynomial algebra in $[\tau]$.
  More precisely, the class 
  \[
    [L(-n_1)\cds L(-n_i)G(-r_1)\cds G(-r_j)\vac ],~~
    n_1\geq \cds \geq n_i\geq 2,~~
    r_1>\cds >r_j\geq \dfr{3}{2},
  \] 
  corresponds to a polynomial in $[\tau]$ of degree at most $2i+j$.
\end{prop}

\pf
We prove that $[L(-n_1)\cds L(-n_i)G(-r_1)\cds G(-r_j)\vac]$ is equivalent to 
a polynomial in $[\tau]$ of degree at most $2i+j$ by induction on the length $i+j$.
By Lemma \ref{lem:a.2} there exists a polynomial $f(X)\in \C[X]$ of degree
at most $i$ such that 
\[
  [L(-n_1)\cds L(-n_i)G(-r_1)\cds G(-r_j)\vac ]
  = [G(-r_1)\cds G(-r_j)\vac ]\,{\ds \mystar_\tw}\,f([\w]).
\]
It follows from Lemma \ref{lem:a.4} that $f([\w])$ is equivalent to a polynomial 
in $[\tau]$ of degree at most $2i$.
By Lemma \ref{lem:a.3} we can rewrite the class $[G(-r_1)\cds G(-r_j)\vac ]$ 
into a sum of shorter monomials in $L(-n)$ and $G(-r)$.
Note that an even element $L(-n)$ appears when we rewrite a pair of two odd elements 
$G(-r)G(-s)$ by taking a commutator.
In the rewriting procedure $L(-n)$ increase the degree of $[\tau]$ in the terminal form 
at most two whereas $G(-r)$ does at most one.
So we can apply the induction and 
$[G(-r_1)\cds G(-r_j)\vac ]$ is equivalent to a polynomial in $[\tau]$ of degree at most $j$.
Therefore, every element of $A_\tw\!\l(\ol{M}_\mathrm{NS}(c,0)\r)$ is equivalent to a polynomial 
in $[\tau]$ of degree as described in the assertion.
This shows that there exists an epimorphism from $\C[X]$ to 
$A_\tw\!\l(\ol{M}_\mathrm{NS}(c,0)\r)$ defined by $X\longmapsto [\tau]$.
We prove that this is the isomorphism.
For any $d \in \C$, the zero-mode $\o(\tau)=G(0)$ has a minimal polynomial $X-d$ on 
the top level of the Verma module $M_{\ram}(c,d)$ over the Ramond algebra.
Hence, $A_\tw\!\l(\ol{M}_\mathrm{NS}(c,0)\r)$ has an irreducible representation 
on which $[\tau]$ acts by an arbitrary scalar.
This implies $A_\tw\!\l(\ol{M}_\mathrm{NS}(c,0)\r)$ is indeed isomorphic to a polynomial 
algebra in $[\tau]$.
\qed
\medskip

Now we describe the twisted Zhu algebra of $L_{\mathrm{NS}}(c_n^1,0)$.

\begin{thm}[\cite{M}]\label{thm:a.6}
  The $\Z_2$-twisted Zhu algebra of $L_{\mathrm{NS}}(c_n^1,0)$ is isomorphic to a 
  quotient of a polynomial ring $\C[X]$ modulo the following polynomial. 
  \[
    \prod_{{1\leq r\leq n+1\atop 1\leq s \leq n+3} \atop r-s\equiv 1 (2)}
    \l( X-\Delta_{r,s}^{(n)}\r). 
  \]
  The isomorphism is given by  $[\tau] \longmapsto X$.
\end{thm}

\pf
By Proposition \ref{prop:a.5}, there is a polynomial $f(X)$ such that 
\[
  A_\tw(L_\mathrm{NS}(c_n^1,0)) \cong \C[X]/\la f(X)\ra .
\]
It follows from the structure of the Verma modules over the Neveu-Schwarz algebra that 
the maximal ideal of $M_\mathrm{NS}(c_n^1,0)$ is generated by two singular vectors, 
$G(-\sfr{1}{2})v_{c_n^1,0}$ and the one of weight $(n+1)(n+3)/2$ (cf.~\cite{IK}).
Let $x$ be the singular vector of $\ol{M}_\mathrm{NS}(c_n^1,0)$ of weight $(n+1)(n+3)/2$.
By Proposition \ref{prop:a.5} there exists a polynomial $g(X)$ such that $g([\tau])=[x]$ in 
$A_\tw\!\l(\ol{M}_\mathrm{NS}(c_n^1,0)\r)$.
It is clear that $f(X)$ divides $g(X)$ since $[x]=g([\tau])=0$ in $A_\tw(L_\ns(c_n^1,0))$.
If $n$ is odd then the weight $(n+1)(n+3)/2$ is an even integer and $x$ is an even element, 
whereas if $n$ is even then $(n+1)(n+3)/2$ is a half-integer and $x$ is an odd element.
By this we see that the possible longest monomial of the form \eqref{eq:a.3} in $x$ is 
$L(-2)^{\fr{1}{4}(n+1)(n+3)}\vac$ if $n$ is odd and $L(-2)^{\fr{1}{4}n(n+4)}G(-\sfr{3}{2})\vac$ 
if $n$ is even, respectively.
Therefore, by Proposition \ref{prop:a.5}, the degree of $g(X)$ is at most $(n+1)(n+3)/2$ 
if $n$ is odd and $n(n+4)/2+1$ if $n$ is even, respectively.
On the other hand, by the GKO construction \cite{GKO} (cf.~Eq.~\eqref{eq:4.7}), 
we have irreducible $\Z_2$-twisted $L_\mathrm{NS}(c_n^1,0)$-modules 
$L_\ram(c_n^1,\Delta_{r,s}^{(n)})$ for $1\leq r\leq n+1$, $1\leq s\leq n+3$ and $r-s\equiv 1$ $(2)$.
The zero-mode $\o(\tau)=G(0)$ acts on the top level of $L_\ram(c_n^1,\Delta_{r,s}^{(n)})$ by 
$\Delta_{r,s}^{(n)}$.
It is straightforward to see that $\Delta_{r,s}^{(n)}$ with $1\leq r\leq n+1$, $1\leq s\leq n+3$, 
$r-s\equiv 1~(2)$ are mutually distinct.
(Note that $h_{r,s,1/2}^{(n)}=h_{n+2-r,n+4-s,1/2}^{(n)}$ but 
$\Delta_{n+2-r,n+4-s}^{(n)}=-\Delta_{r,s}^{(n)}$.)
Therefore, $f(X)$ is divisible by the polynomial
\[
  \prod_{{1\leq r\leq n+1\atop 1\leq s \leq n+3} \atop r-s\equiv 1 (2)} (X-\Delta_{r,s}^{(n)}).
\]
The degree of the polynomial above is $(n+1)(n+3)/2$ if $n$ is odd and 
$n(n+4)/2+1$ if $n$ is even.
Therefore, by comparing degrees, we see that both $g(X)$ and $f(X)$ are scalar multiples 
of the polynomial above.
This completes the proof.
\qed
\medskip

As a corollary, we obtain the classification of irreducible $\Z_2$-twisted $L_\ns(c_n^1,0)$-modules.

\begin{thm}[\cite{M}]\label{thm:a.7}
  The irreducible $\Z_2$-twisted $L_{\mathrm{NS}}(c_n^1,0)$-modules are 
  $L_\ram(c_n^1,\Delta_{r,s}^{(n)})$, $1\leq r\leq n+1$, $1\leq s\leq n+3$, $r-s\equiv 1 \mod 2$.
\end{thm}

\begin{thm}\label{thm:a.8}
  $L_{\mathrm{NS}}(c_n^1,0)$ is $\Z_2$-rational.
\end{thm}

\pf
In this proof we use the notation as in Section \ref{sec:2.2}.
Since the even part $L_\ns(c_n^1,0)^{[0]}$ is $C_2$-cofinite by (2) of Theorem \ref{thm:4.3}, 
every $\Z_2$-twisted $L_\ns(c_n^1,0)$-module is $\N$-gradable by \cite{ABD,Mi2}.
Let $M=\oplus_{n\geq 0}M(n)$ be an $\N$-graded $\Z_2$-twisted $L_\ns(c_n^1,0)$-module with 
non-trivial top level $M(0)$.
Then by Theorems \ref{thm:a.1} and \ref{thm:a.6} $M(0)$ is a semisimple 
$A_\tw(L_\ns(c_n^1,0))$-module and is a direct sum of eigenvectors of 
$\o(G(-\sfr{3}{2})\vac)=G(0)$ with eigenvalues 
$\Delta_{r,s}^{(n)}$, $1\leq r\leq n+1$, $1\leq s\leq n+3$, $r-s\equiv 1 \mod 2$.
We shall show that every eigenvector of $G(0)$ generates an irreducible $\Z_2$-twisted submodule.
Let $x$ be a $G(0)$-eigenvector of $M$ with eigenvalue $\Delta_{r,s}^{(n)}$ and 
let $X$ be the submodule generated by $x$.
Then up to linearity there is a unique epimorphism 
$\pi: M_\ram(c_n^1,\Delta_{r,s}^{(n)})\longrightarrow X$.
It is shown in Theorem 4.2 of \cite{IK} that every submodule of $M_\ram(c_n^1,\Delta_{r,s}^{(n)})$ 
is generated by singular vectors\footnote{%
In the Ramond case the Verma module of central charge $c$ and highest weight $h$ in \cite{IK} 
corresponds to 
$M_\ram(c,\sqrt{h-c/24})\oplus M_\ram(c,-\sqrt{h-c/24})$ if $h \ne c/24$, 
and to an extension of $M_\ram(c,0)$ by itself if $h=c/24$.}.
For $1\leq r\leq n+1$ and $1\leq s\leq n+3$, set
\[
  h_{r,s,1/2}^{(n)}(i)
  :=\begin{cases}
    \dfr{((i(n+2)+r)(n+4)-s(n+2))^2-4}{8(n+2)(n+4)}+\dfr{1}{16} 
    & \mbox{ if }~ i\equiv 0 \mod 2,
    \medskip\\
    \dfr{(((i-1)(n+2)+r)(n+4)+s(n+2))^2-4}{8(n+2)(n+4)}+\dfr{1}{16} 
    & \mbox{ if }~ i\equiv 1 \mod 2.
  \end{cases}
\]
It is also shown in (loc.~cit.) that the $L(0)$-weights of singular vectors of 
$M_\ram(c_n^1,\Delta_{r,s})$ belong to the set 
\[
  H_{r,s,1/2}^{(n)}=\l\{ h_{r,s,1/2}^{(n)}(i) \,\Big|~ i\in \Z\setminus \{ 0\}\r\} .
\]
Suppose $X$ is reducible and we have an $L(0)$-homogeneous singular vector $y$ of $X$.
Then by Theorem \ref{thm:a.6} the $L(0)$-weight of $y$ is equal to 
$h_{r',s',1/2}^{(n)}=(\Delta_{r',s'}^{(n)})^2+c_n^1/24$ for some $1\leq r'\leq n+1$, 
$1\leq s'\leq n+3$ and $r'-s'\equiv 1 \mod 2$.
However, for $1\leq r,r'\leq n+1$ and $1\leq s,s'\leq n+3$, it is directly verified that 
$h_{r',s',1/2}^{(n)}\not\in H_{r,s,1/2}^{(n)}$.
Thus $X$ has no singular vector and is irreducible.
Now the theorem follows from Proposition 5.11 of \cite{DLTYY}.
\qed


\begin{thebibliography}{AAAAA12}

\bibitem[ABD04]{ABD}
  T. Abe, G. Buhl and C. Dong, 
  Rationality, regularity, and $C_2$-cofiniteness.
  \textit{Trans. Amer. Math. Soc.} \textbf{356} (2004), 3391--3402.

\bibitem[A97]{A1}
  D. Adamovi\'{c}, 
   Rationality of Neveu-Schwarz vertex operator superalgebras. 
  \textit{Internat. Math. Res. Notices} \textbf{17}, 865--874.

\bibitem[A04]{A2}
  D. Adamovi\'{c}, Regularity of certain vertex operator superalgebras.
  \textit{Contemp. Math.} \textbf{343} (2004), 1--16. 

\bibitem[DJX13]{DJX}
  C. Dong, X. Jiao and F. Xu, 
  Quantum dimensions and quantum Galois theory. 
  \textit{Trans. Amer. Math. Soc.} \textbf{365} (2013), no. 12, 6441--6469. 

\bibitem[DLMN98]{DLMN}
  C. Dong, H. Li, G. Mason and S.P. Norton, 
 Associative subalgebras of Griess algebra and related topics.
 Proc. of the Conference on the Monster and Lie algebra at the Ohio State
  University, May 1996, ed. by J. Ferrar and K. Harada, Walter de Gruyter, 
  Berlin - New York, 1998.

\bibitem[DLTYY04]{DLTYY}
  C. Dong, C.H. Lam, K. Tanabe, H. Yamada and K. Yokoyama, 
  $\Z_3$ symmetry and $W_3$ algebra in lattice vertex operator algebras.
  \textit{Pacific J. Math.} \textbf{215} (2004), 245--296.

\bibitem[DMZ94]{DMZ}
  C. Dong, G. Mason and Y. Zhu, 
  Discrete series of the Virasoro algebra and the moonshine module.
  Proc. Symp. Pure. Math., American Math. Soc. \textbf{56} II (1994), 295--316.

\bibitem[DZ06]{DZ}
  C.Dong and Z. Zhao,  
  Twisted representations of vertex operator superalgebras. 
  \textit{Commun. Contemp. Math.} \textbf{8} (2006), no. 1, 101--121. 

\bibitem[FZ92]{FZ}
  I.B. Frenkel and Y. Zhu, 
  Vertex operator algebras associated to representation of affine and Virasoro algebras. 
  \textit{Duke Math. J.} \textbf{66} (1992), 123--168.

\bibitem[GKO86]{GKO}
  P. Goddard, A. Kent and D. Olive, 
  Unitary representations of the Virasoro and super-Virasoro algebras. 
  \textit{Comm. Math. Phys.} \textbf{103} (1986), 105--119.

\bibitem[IK03]{IK}
  K. Iohara and Y. Koga, 
  Representation theory of Neveu-Schwarz and Ramond algebras I: Verma modules.
  \textit{Adv. Math.} \textbf{178} (2003), 1--65.

\bibitem[HLY12a]{HLY1}
  G. Hoehn, C.H. Lam and H. Yamauchi, 
  McKay's $E_7$ observation on the Babymonster. 
  \textit{Internat. Math. Res. Notices} \textbf{2012} (2012), 166--212.

\bibitem[HLY12b]{HLY2}
  G. Hoehn, C.H. Lam and H. Yamauchi,
  McKay's $E_6$ observation on the largest Fischer group. 
  \textit{Comm. Math. Phys.} \textbf{310} Vol. 2 (2012), 329--365. 

\bibitem[K98]{K}
  V. Kac, 
  Vertex algebras for beginners, second edition. 
  Cambridge University Press, Cambridge, 1998.

\bibitem[KR86]{KR}
  V. Kac and A.K. Reina, 
  Bombay lectures on highest weight representations of infinite-dimensional Lie algebras. 
  Advanced Series in Mathematical Physics, 2. (1986), World Scientific Publishing.

\bibitem[KW94]{KW}
  V. Kac and W. Wang, 
  Vertex operator superalgebras and their representations.
  Mathematical Aspects of Conformal and Topological Field Theories and Quantum Groups
  (South Hadley, MA, 1992), Contemp. Math. 175, Amer. Math. Soc., Providence, 1994, 
  161--191.
 \bibitem[LS08]{LS}
  C.H. Lam and S. Sakuma, 
  On a class of vertex operator algebras having a faithful $S_{n+1}$-action. 
  \textit{Taiwanese J. Math.} \textbf{12} (2008), no. 9, 2465--2488. 
   
  \bibitem[LYY03]{LYY1}
  C.H. Lam, H. Yamada and H. Yamauchi, Vertex operator algebras, extended
  $E_8$-diagram, and McKay's observation on the Monster simple group.
  {\it Trans. Amer. Math. Soc.} {\bf 359} (2007), 4107--4123.

\bibitem[LYY05]{LYY2}
  C.H. Lam, H. Yamada and H. Yamauchi, McKay's observation and vertex
  operator algebras generated by two conformal vectors of central
  charge 1/2. {\it Internat. Math. Res. Papers} {\bf 3} (2005), 117--181.

\bibitem[LY16]{LY}
  C.H. Lam and H. Yamauchi, 
  The Conway-Miyamoto correspondences for the Fischer 3-transposition groups.
  Preprint.


\bibitem[Li96]{Li1}
  H. Li, 
  Local systems of twisted vertex operators, vertex operator superalgebras and twisted modules.  
  \textit{J. Pure Appl. Algebra} \textbf{109} (1996), no. 2, 143--195. 

\bibitem[Li97]{Li2}
  H. Li, 
  Extension of vertex operator algebras by a self-dual simple module.
  \textit{J. Algebra} \textbf{187} (1997), no. 1, 236--267. 

\bibitem[M07]{M}
  A. Milas, 
  Characters, supercharacters and Weber modular functions. 
  \textit{J. Reine Angew. Math.} \textbf{608} (2007) 35--64. 

\bibitem[Mi96]{Mi1}
  M. Miyamoto, 
  Griess algebras and conformal vectors in vertex operator algebras. 
  \textit{J. Algebra} \textbf{179} (1996), 528--548.
  
\bibitem[Mi04]{Mi4}
  M. Miyamoto, A new construction of the moonshine vertex operator
  algebra over the real number field. {\it Ann. of Math.} {\bf 159}
  (2004), 535--596.
  
\bibitem[Mi04b]{Mi2}
  M. Miyamoto,
  Modular invariance of vertex operator algebras satisfying $C_2$-cofiniteness.
  \textit{Duke. Math. J.} \textbf{122} (2004), 51--91.
  
  


\bibitem[Mi14]{Mi3}
  M. Miyamoto, 
  $C_2$-cofiniteness of cyclic-orbifold model.
  \textit{Comm. Math. Phys.} \textbf{335} (2015), 1279--1286.
  
  
\bibitem[Sa07]{Sak}
  S. Sakuma, 6-transposition property of $\tau$-involutions of
  vertex operator algebras.
  {\it Internat. Math. Res. Notices} {\bf 2007},  no. {\bf 9},
  Art. ID rnm 030, 19 pp.

\bibitem[W93]{W}
  W. Wang, 
  Rationality of Virasoro vertex operator algebras.
  \textit{Internat. Math. Res. Notices} \textbf{71} (1993), 197--211.

\bibitem[Y05]{Y1}
  H. Yamauchi, 
  2A-orbifold construction and the baby-monster  vertex operator superalgebra. 
  \textit{J. Algebra} \textbf{284} (2005), 645--668.


\bibitem[Z96]{Z}
  Y. Zhu, 
  Modular invariance of characters of vertex operator algebras. 
  \textit{J. Amer. Math. Soc.} \textbf{9} (1996), 237--302.


\end{thebibliography}
\end{document}